\documentclass[12pt, 3p
			  ]{elsarticle}
\usepackage{geometry}
\usepackage{amsmath,amsbsy,amsfonts}
\allowdisplaybreaks[4]
\usepackage{bm}
\usepackage{amssymb}

\usepackage{booktabs}	
\usepackage[colorlinks]{hyperref}
\usepackage{enumerate}
\usepackage{lineno}
\usepackage{mathrsfs}
\usepackage{stmaryrd} % \jump
\usepackage{latexsym}
\usepackage{caption,subcaption}
\usepackage{graphicx}
\usepackage{graphics}
\usepackage{epstopdf}
\usepackage{ntheorem}
\usepackage{pstricks-add}
\usepackage{CJKutf8,CJKnumb}%,CJKulem}
\AtBeginDocument{\begin{CJK}{UTF8}{gbsn}}
\AtEndDocument{\clearpage \end{CJK}}
\usepackage{cases}
\usepackage{float}
\usepackage{diagbox}
\newfam\msbfam
\font\tenmsb=msbm10    \textfont\msbfam=\tenmsb
\font\sevenmsb=msbm7 \scriptfont\msbfam=\sevenmsb
\font\fivemsb=msbm5 \scriptscriptfont\msbfam=\fivemsb

\newfam\bigfam
\font\tenbig=msbm10 scaled \magstep2   \textfont\bigfam=\tenbig
\font\sevenbig=msbm7 scaled \magstep2 \scriptfont\bigfam=\sevenbig
\font\fivebig=msbm5 scaled \magstep2
\scriptscriptfont\bigfam=\fivebig

\definecolor{mygray}{gray}{.85}
\numberwithin{equation}{section}

\journal{}

\begin{document}

\newtheorem{theorem}{Theorem}[section]
\newtheorem{definition}{Definition}[section]
\newtheorem{example}{Example}[section]
\newtheorem{lemma}[theorem]{Lemma}
\newtheorem{corollary}[theorem]{Corollary}
\newtheorem{proposition}[theorem]{Proposition}
\newtheorem{remark}{Remark}[section]
\theorembodyfont{\normalfont}
\newproof{pf}{Proof}
\newproof{pot}{Proof of Theorem \ref{thm2}}
\def\endpf{\hspace*{\fill}~$\square$\par\endtrivlist\unskip}

\begin{frontmatter}

\title{On the solution of the coupled steady-state dual-porosity-Navier-Stokes fluid flow model with the Beavers-Joseph-Saffman interface condition\tnoteref{mytitlenote}}
%\tnotetext[mytitlenote]{This work is partially supported by NSF grant DMS-1720212,
%National Natural Science Foundation of China(11371289, 11471329, 11771348), National Magnetic Confinement Fusion Science Program of China (2015GB110003), National Key Research and Development Program of China (2016YFB0201304),  Major Research and Development Program of China (2016YFB0200901) and Youth Innovation Promotion Association of CAS (2016003).}

\author[add1]{Di Yang}\ead{yangdi0226@163.com}
\author[add1]{Yinnian He\corref{correspondingauthor}}\ead{heyn@mail.xjtu.edu.cn}
%\author[add1]{Shuaichao Pei}\ead{peishuaichao215@stu.xjtu.edu.cn}
\author[add1]{Luling Cao}\ead{lulingcao@163.com}
\cortext[correspondingauthor]{Corresponding author.}

\address[add1]{School of Mathematics and Statistics, Xi'an Jiaotong University, Xi'an, Shaanxi 710049, P.R. China }

\begin{abstract}
\par In this work, we propose a new analysis strategy to establish an a priori estimate of the weak solutions to the coupled steady-state dual-porosity-Navier-Stokes fluid flow model with the Beavers-Joseph-Saffman interface condition. The most advantage of our proposed method is that the a priori estimate and the existence result are independent of small data and the large viscosity restriction. Therefore the global uniqueness of the weak solution is naturally obtained.

\end{abstract}

\end{frontmatter}

% \tableofcontents
%\linenumbers

\section{Introduction}
\label{section:introuction}

Coupled free flow and porous medium flow systems play an important role in many practical engineering fields, e.g., the flood simulation of arid areas in geological science~\cite{Discacciati2004phd}, filtration treatment in industrial production~\cite{ hanspal2006numerical, nassehi1998modelling}, petroleum exploitation in mining, and blood penetration between vessels and organs in life science~\cite{d2011robust}. Specifically, the systems are usually described by Navier-Stokes equations (or Stokes equations) coupled with Darcy's equation, and there are amounts of achievements such as \cite{Cao20132013, Cao2011Robin, Cai2009NUMERICAL, Chen2012Efficient, Girault20092052, Mu2010Decoupled, zhao2016two-grid}. However, the standard Darcy's equation describes fluids flowing through only a single porosity medium, which is not accurate to deal with the complicated multiple porous media similar to naturally fractured reservoir. Actually, the naturally fractured reservoir is comprised of low permeable rock matrix blocks surrounded by an irregular network of natural microfractures, and further they have different fluid storage and conductivity properties~\cite{New1983, The1963}. In 2016, Hou et al. \cite{Hou2016710} proposed and numerically solved a coupled dual-porosity-Stokes fluid flow model with four multi-physics interface conditions. The authors used the dual-porosity equations over Darcy's region to describe fluid flowing through the multiple porous medium. Recently, several related research on the above model can be found in the literatures~\cite{AlMahbub2019803, AlMahbub2020112616, Gao202001, he2020an, Shan2019389}. In particular, Gao and Li \cite{Gao202001} proposed a decoupled stabilized finite element method to solve the coupled dual-porosity-Navier-Stokes fluid flow model in the numerical field.

The steady-state dual-porosity-Navier-Stokes fluid flow model has distinct features and difficulties in mathematical analysis. Many numerical methods have been studied for the well-known stationary or time-dependent Navier-Stokes/Darcy model with Beavers-Joseph or Beaver-Joseph-Saffman interface condition, including coupled finite element methods \cite{Badea2010195, Cao2021113128, Discacciati2009315, Discacciati2017571, Zuo20151009}, discontinuous Galerkin methods \cite{Chidyagwai20093806, Girault20092052, Girault201328}, domain decomposition methods \cite{Cao20132013, Discacciati2004phd, Discacciati2009315, He2015264, Qiu2020109400}, and decoupled methods based on two-grid finite element \cite{Cai2009NUMERICAL, Fang2020915, zhao2016two-grid, Zuo20151009}. In spite of the above great contributions to numerical simulation, the existence of a weak solution to the coupled dual-porosity-Navier-Stokes fluid flow model with Beavers-Joseph-Saffman interface condition for general data keeps unresolved. In many literatures \cite{Badea2010195, Cai2009NUMERICAL, Chidyagwai20093806, Discacciati2004phd, Discacciati2009315, Discacciati2017571, Girault20092052, Girault201328, He2015264, Zuo20151009}, a priori estimates and existence of a weak solution need suitable small data and/or large viscosity restrictions, and therefore only local uniqueness can be established when the data satisfy additional restrictions. In \cite{Girault20092052}, the authors pointed out that the difficulty for a priori estimates and existence with general data is stemmed from the transmission interface condition, which does not completely compensate the nonlinear convection term from the Navier-Stokes equations in the energy balance.

Therefore in this paper, stemming from resolving steady-state Navier-Stokes equations with mixed boundary conditions in \cite{Hou201947}, we shall establish a new a priori estimate of the weak solutions by coupling the model problem with a designed auxiliary problem in order to completely compensate the nonlinear convection term from the Navier-Stokes equations. In addition, we shall also prove existence of a weak solution without small data or large viscosity restriction. As a result, the global uniqueness of the weak solution is naturally obtained.

The rest of this paper is organized as follows. In Section 2, we specify the steady-state dual-porosity-Navier-Stokes fluid flow model with Beavers-Joseph-Saffman interface condition and provide its Galerkin variational formulation. In Section 3, we establish a new a priori estimate of the weak solutions by coupling the model problem with an auxiliary problem, which is designed subject to the model problem. Finally, in Section 4, we prove existence of the weak solution without small data and/or large viscosity restriction by the Galerkin method and Brouwer's fixed-point theorem, and global uniqueness of all variables by the inf-sup condition and Babu\u{s}ka--Brezzi's theory.

\section{Model specification}
\label{section:model}

\subsection{Setting of the problem}
In this section, we consider the steady-state dual-porosity-Navier-Stokes fluid flow model in a bounded open polygonal domain $\Omega\subset\mathbb{R}^N(N=2,3)$ with four physically valid interface conditions.

\begin{figure}[!ht]
  \centering
  \includegraphics[width=0.9\textwidth]{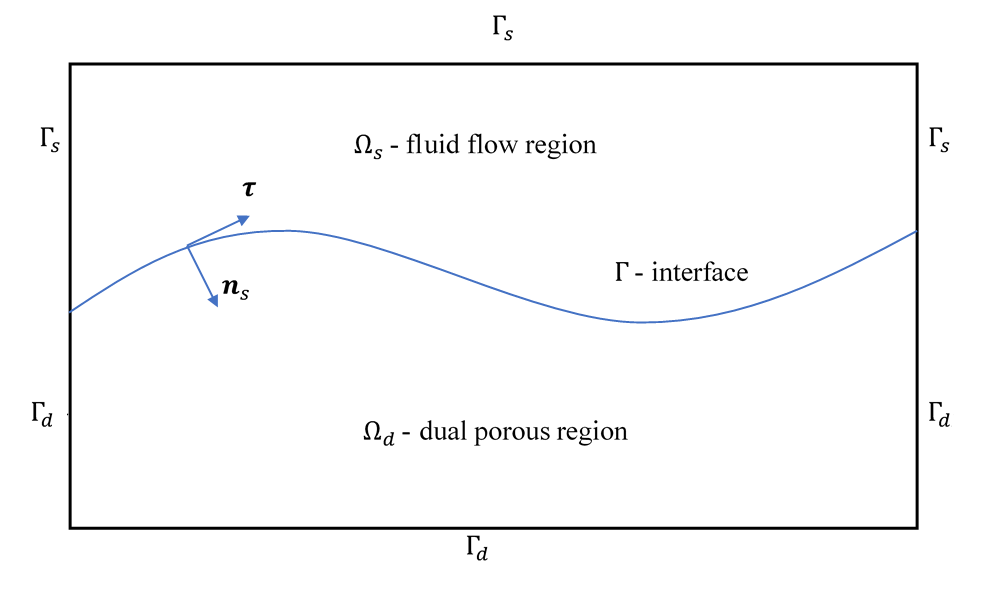}\\
  \caption{Schematic representation of fluid flow region $\Omega_s$ and dual-porous medium region $\Omega_d$, separated by an interface $\Gamma$.}\label{global:domain}
\end{figure}

The domain $\Omega$ consists of a fluid flow region  $\Omega_s$ and a dual-porous medium region $\Omega_d$, with interface $\Gamma=\overline{\Omega}_s\cap\overline{\Omega}_d$ (see Figure \ref{global:domain}). Both $\Omega_s$ and $\Omega_d$ are open, regular, simply connected, and bounded with Lipschitz continuous boundaries $\Gamma_i=\partial\Omega_i\setminus\Gamma~ (i=s,d)$, respectively. Here, $\Omega_s\cap\Omega_d=\emptyset$, $\overline{\Omega}_s\cup\overline{\Omega}_d=\Omega$. The unit normal vector of the interface $\Gamma$ pointing from $\Omega_s$ to $\Omega_d$ (from $\Omega_d$ to $\Omega_s$) is denoted by $\bm{n}_s$ (resp. $\bm{n}_d$), and the corresponding unit tangential vectors are denoted by $\bm{\tau}_i$, $i=1,\cdots,N-1$. In two-dimensional case, if we write $\bm{n}_s=(n_1,n_2)\in\mathbb{R}^2$, then $\bm{\tau}=\bm{n}_s^\top=(-n_2,n_1)$.

In $\Omega_s$, the fluid flow is governed by the Navier-Stokes equation \cite{Cao2021113128, Cao20132013, Chidyagwai20093806, Discacciati2004phd, Discacciati2009315, Girault20092052, Girault201328}:

\begin{equation}\label{stokes}
\begin{cases}
\left(\bm{u}_s\cdot\nabla\right)\bm{u}_s-\nabla\cdot
\mathbb{T}(\bm{u}_s,p_s)=\bm{f}_s&\text{in}~\Omega_s,\\
\nabla\cdot\bm{u}_s=0&\text{in}~\Omega_s,\\
\bm{u}_s=\bm{0}&\text{on}~\Gamma_s,
\end{cases}
\end{equation}
where $\mathbb{T}(\bm{u}_s,p_s)=-p_s\mathbb{I}+2\nu\mathbb{D}(\bm{u}_s)$ is the stress tensor, $\nu$ is the kinetic viscosity, $\mathbb{D}(\bm{u}_s)=\frac{1}{2}\left[\nabla\bm{u}_s+(\nabla\bm{u}_s)^\top\right]$ is the velocity deformation tensor, $\bm{u}_s(\bm{x})$ denotes the fluid velocity in $\Omega_s$, $p_s(\bm{x})$ denotes the kinematic pressure in $\Omega_s$, and $\bm{f}_s(\bm{x})$ denotes a general body force term that includes gravitational acceleration. As usual, we write formally:

\[(\bm{v}\cdot\nabla)\bm{w}=\sum_{i=1}^N v_i\frac{\partial\bm{w}}{\partial x_i},\quad
\nabla\cdot\bm{v}=\sum_{i=1}^N\frac{\partial v_i}{\partial x_i}.\]

The filtration of an incompressible fluid through porous media is often described using Darcy's law. So in dual-porous medium domain $\Omega_d$, the flow is governed by a traditional dual-porosity model, which is composed of matrix and microfracture equations as follows:

\begin{equation}\label{dual:darcy}
\begin{cases}
-\nabla\cdot
\big(\frac{k_m}{\mu}\nabla\phi_m\big)+\frac{\sigma k_m}{\mu}(\phi_m-\phi_f)=0&\text{in}~\Omega_d,\\
-\nabla\cdot
\big(\frac{k_f}{\mu}\nabla\phi_f\big)+\frac{\sigma k_m}{\mu}(\phi_f-\phi_m)=f_d&\text{in}~\Omega_d,\\
\phi_m=0&\text{on}~\Gamma_{d},\\
\phi_f=0&\text{on}~\Gamma_{d},
\end{cases}
\end{equation}
where $\phi_m(\bm{x})$, $\phi_f(\bm{x})$ denote the matrix and microfracture flow pressure, $\mu$ is the dynamic viscosity, $k_m$, $k_f$ are the intrinsic permeability of the matrix and microfracture regions, $\sigma$ is the shape factor characterizing the morphology and dimension of the microfractures, $f_d$ is the sink/source term for the microfractures, and the term $\frac{\sigma k_m}{\mu}(\phi_m-\phi_f)$ describes the mass exchange between matrix and microfractures.

%\begin{remark}
%We impose homogeneous Dirichlet boundary conditions for $\bm{u}$, $\phi_m$ and $\phi_f$ in order to analyse. However in practical applications, more realistic boundary conditions may need to be imposed, which will be discussed in Section .
%\end{remark}

Based on the fundamental properties of the dual-porosity fluid flow model and traditional Stokes--Darcy flow model, Hou et al. \cite{Hou2016710} introduced four physically valid interface conditions as follows to couple appropriately the dual-porosity-Stokes model, which are also adapted to our problem.

\begin{itemize}
  \item No-exchange condition between the matrix and the conduits/macrofractures:
  \begin{equation}\label{interface:condition:1}
  -\frac{k_m}{\mu}\nabla\phi_m\cdot\bm{n}_d=0
  \quad\text{on}~\Gamma.
  \end{equation}

  \item Mass conservation:
  \begin{equation}\label{interface:condition:2}
  \bm{u}_s\cdot\bm{n}_s=\frac{k_f}{\mu}
  \nabla\phi_f\cdot\bm{n}_d
  \quad\text{on}~\Gamma.
  \end{equation}

  \item Balance of normal forces:
  \begin{equation}\label{interface:condition:3}
  -\bm{n}_s\cdot(\mathbb{T}(\bm{u}_s,p_s)\bm{n}_s)
  =\frac{\phi_f}{\rho}
  \quad\text{on}~\Gamma.
  \end{equation}

  \item The Beavers-Joseph-Saffman interface condition \cite{Beavers1967197, Jone1973231, Mikelic20001111, Saffman197193}: for $i=1,\cdots,N-1$,
  \begin{equation}\label{interface:condition:4}
  -\bm{\tau}_i\cdot\left(\mathbb{T}(\bm{u}_s,p_s)\bm{n}_s\right)
  =\frac{\alpha\nu\sqrt{N}}{\sqrt{\text{trace}(\bm{\Pi})}}
  \bm{u}_s\cdot\bm{\tau}_i
  \quad \text{on}~\Gamma.
  \end{equation}
\end{itemize}
In \eqref{interface:condition:1}--\eqref{interface:condition:4}, $\alpha$ is the Beavers-Joseph constant depending on the properties of the dual-porous medium, $\bm{\Pi}$ represents the intrinsic permeability that satisfies the relation $\bm{\Pi}=k_f\mathbb{I}$, $\mathbb{I}$ is the $N\times N$ identity matrix, and $\rho$ is the fluid density.

\subsection{Galerkin variational formulation}

Throughout this paper we use the following standard function spaces. For a Lipschitz domain $\mathcal D\subset\mathbb{R}^N$, $N\geqslant 1$, we denote by $W^{k,p}(\mathcal D)$ the Sobolev space with indexes $k\geqslant 0$, $1\leqslant p\leqslant\infty$ of real-valued functions defined on $\mathcal D$, endowed with the seminorm $|\cdot|_{W^{k,p}(\mathcal D)}$ denoted by $|\cdot|_{k,p,\mathcal D}$ and norm $\|\cdot\|_{W^{k,p}(\mathcal D)}$ denoted by $\|\cdot\|_{k,p,\mathcal D}$ \cite{Adams2003}. When $p=2$, $W^{s,2}(\mathcal D)$ is denoted as $H^k(\mathcal D)$ and the corresponding seminorm and norm are written as $|\cdot|_{k,\mathcal D}$ and $\|\cdot\|_{k,\mathcal D}$, respectively. In addition, with $|\mathcal D|$ we denote the $N$-dimensional Hausdorff measure of $\mathcal D$.

To perform the variational formulation, we define some necessary Hilbert spaces given by

\begin{align*}
\bm{X}_s&:=\left\{
\bm{v}\in\bm{H}^1(\Omega_s):\bm{v}=\bm{0}~\text{on}~\Gamma_s
\right\},\\
X_d&:=\left\{
\psi\in H^1(\Omega_d):\psi=0~\text{on}~\Gamma_d
\right\},\\
Q_s&:=L^2(\Omega_s).
\end{align*}
We also need the trace space $\bm{H}_{00}^{1/2}(\Gamma):=\bm{X}_s|_\Gamma$ (resp. $H_{00}^{1/2}(\Gamma):=X_d|_\Gamma$), which is a nonclosed subspace of $\bm{H}^{1/2}(\Gamma)$ (resp. $H^{1/2}(\Gamma)$) and has a continuous zero extension to $\bm{H}^{1/2}(\partial\Omega_s)$ (resp. $H^{1/2}(\partial\Omega_d)$) \cite{Cao20104239, Cao20101, Discacciati2004phd}. For the trace space $\bm{H}_{00}^{1/2}(\Gamma)$ and its dual space $(\bm{H}_{00}^{1/2}(\Gamma))'$, we have the following continuous imbedding result \cite{Cao20101}:
\begin{equation}\label{interface:space:imbedding}
\bm{H}_{00}^{1/2}(\Gamma)
\subsetneqq \bm{H}^{1/2}(\Gamma)
\subsetneqq \bm{L}^2(\Gamma)
\subsetneqq \bm{H}^{-1/2}(\Gamma)
\subsetneqq (\bm{H}_{00}^{1/2}(\Gamma))'.
\end{equation}
One can see more details in \cite{Cao20101, Li201692, Shan2013813} and the references therein. For any bounded domain $\mathcal D\subset\mathbb{R}^N$, $(\cdot,\cdot)_{\mathcal D}$ denotes the $L^2$ inner product on $\mathcal D$, and $\langle\cdot,\cdot\rangle_{\partial\mathcal D}$ denotes the $L^2$ inner product (or duality pairing) on the boundary $\partial\mathcal D$. We also consider the following product Hilbert space
$\underline{\bm{Y}}:=\bm{X}_s\times X_d\times X_d$ with norm

\[
\|\underline{\bm{w}}\|_{\underline{\bm{Y}}}=\Big(\|\bm{w}_s\|_{1,\Omega_s}^2+\|\psi_m\|_{1,\Omega_d}^2
+\|\psi_f\|_{1,\Omega_d}^2\Big)^{1/2},\quad\forall\,\underline{\bm{w}}=(\bm{w}_s,\psi_m,\psi_f)\in\underline{\bm{Y}}.
\]
In addition, based on the following formula:

\[
\left((\bm{u}\cdot\nabla)\bm{v},\bm{w}\right)_{\mathcal D}
=\left\langle\bm{u}\cdot\bm{n},\bm{v}\cdot\bm{w}\right\rangle_{\partial\mathcal D}
-\left((\bm{u}\cdot\nabla)\bm{w},\bm{v}\right)_{\mathcal D}-\left((\nabla\cdot\bm{u})\bm{v},\bm{w}\right)_{\mathcal D},\quad\forall\,\bm{u},\bm{v},\bm{w}\in\bm{H}^1(\mathcal D),
\]
we introduce the trilinear form $b(\cdot;\cdot,\cdot)$ given by $\forall\,\bm{u}_s,\bm{v}_s,\bm{w}_s\in\bm{X}_s$,

\begin{equation}\label{trilinear:form}
\begin{split}
b(\bm{u}_s;\bm{v}_s,\bm{w}_s)
&=\left((\bm{u}_s\cdot\nabla)\bm{v}_s,\bm{w}_s\right)_{\Omega_s}
+\frac{1}{2}\left((\nabla\cdot\bm{u}_s)\bm{v}_s,\bm{w}_s\right)_{\Omega_s}\\
&=\frac{1}{2}\left\langle\bm{u}_s\cdot\bm{n}_s,\bm{v}_s\cdot\bm{w}_s\right\rangle_\Gamma
+\frac{1}{2}\left((\bm{u}_s\cdot\nabla)\bm{v}_s,\bm{w}_s\right)_{\Omega_s}
-\frac{1}{2}\left((\bm{u}_s\cdot\nabla)\bm{w}_s,\bm{v}_s\right)_{\Omega_s}.
\end{split}
\end{equation}

Hence, the Galerkin variational formulation of the coupled problem \eqref{stokes}--\eqref{interface:condition:4} is proposed that: to find $(\underline{\bm{u}},p_s)\in\underline{\bm{Y}}\times Q_s$ such that

\begin{equation}\label{continuous:weak:formulation}
a(\underline{\bm{u}},\underline{\bm{v}})
+d(\bm{v}_s,p_s)
-d(\bm{u}_s,q_s)
+b(\bm{u}_s;\bm{u}_s,\bm{v}_s)
=\rho(\bm{f}_s,\bm{v}_s)_{\Omega_s}
+(f_d,\psi_f)_{\Omega_d},\quad
\forall\,(\underline{\bm{v}},q_s)\in\underline{\bm{Y}}\times Q_s,
\end{equation}
where $\underline{\bm{u}}=(\bm{u}_s,\phi_m,\phi_f)$, $\underline{\bm{v}}=(\bm{v}_s,\psi_m,\psi_f)$, the bilinear forms $a(\cdot,\cdot)$ and $d(\cdot,\cdot)$ are defined as

\begin{align*}
a_s(\underline{\bm{u}},\underline{\bm{v}})&=
2\rho\nu\left(\mathbb{D}(\bm{u}_s),\mathbb{D}(\bm{v}_s)\right)_{\Omega_s},\\
a_d(\underline{\bm{u}},\underline{\bm{v}})&=
\frac{k_m}{\mu}\left(\nabla\phi_m,\nabla\psi_m\right)_{\Omega_d}+
\frac{k_f}{\mu}\left(\nabla\phi_f,\nabla\psi_f\right)_{\Omega_d}+
\frac{\sigma k_m}{\mu}\left(\phi_m-\phi_f,\psi_m\right)_{\Omega_d}+
\frac{\sigma k_m}{\mu}\left(\phi_f-\phi_m,\psi_f\right)_{\Omega_d},\\
a_\Gamma(\underline{\bm{u}},\underline{\bm{v}})&=
\left\langle\phi_f,\bm{v}_s\cdot\bm{n}_s\right\rangle_\Gamma
-\left\langle\psi_f,\bm{u}_s\cdot\bm{n}_s\right\rangle_\Gamma
+\sum_{i=1}^{N-1}\frac{\alpha\rho\nu}{\sqrt{k_f}}
\left\langle\bm{u}_s\cdot\bm{\tau}_i,\bm{v}_s\cdot\bm{\tau}_i\right\rangle_\Gamma,\\
a(\underline{\bm{u}},\underline{\bm{v}})&=a_s(\underline{\bm{u}},\underline{\bm{v}})+
a_d(\underline{\bm{u}},\underline{\bm{v}})+a_\Gamma(\underline{\bm{u}},\underline{\bm{v}}),\\
d(\bm{v}_s,q_s)&=-\rho(\nabla\cdot\bm{v}_s,q_s)_{\Omega_s}.
\end{align*}

\begin{remark}
We note that the term $\frac{1}{2}\left((\nabla\cdot\bm{u}_s)\bm{v}_s,\bm{w}_s\right)_{\Omega_s}$ vanishes in \eqref{trilinear:form} if $\nabla\cdot\bm{u}_s=0$.
\end{remark}

\begin{remark}\label{remark:dirichlet:boundary}
For the sake of clarity, in our analysis we shall adopt homogenous boundary conditions. In addition, for the general Dirichlet boundary conditions:

\begin{equation*}%\label{general:Dir:bc}
\bm{u}|_{\Gamma_s}=\bm{u}^{\text{dir}},\quad
\phi_m|_{\Gamma_d}=\phi_m^{\text{dir}},\quad
\phi_f|_{\Gamma_d}=\phi_f^{\text{dir}},
\end{equation*}
the standard homogenization technique and the lifting operators employed in \cite{Discacciati2004phd} can be used to obtain an equivalent system with the homogeneous Dirichlet boundary conditions.
\end{remark}

Thanks to \cite{Chidyagwai20093806, Girault20092052}, we have the following important result:
\begin{lemma}\label{lemma:solution:equvi}
Assume that $\bm{f}_s\in\bm{L}^2(\Omega_s)$ and $f_d\in L^2(\Omega_d)$. Then if $(\bm{u}_s,\phi_f,\phi_m,p_s)\in \bm{X}_s\times X_d\times X_d\times Q_s$ satisfies \eqref{stokes}--\eqref{interface:condition:4}, it is also a solution to problem \eqref{continuous:weak:formulation}. Conversely, any solution of problem \eqref{continuous:weak:formulation} satisfies \eqref{stokes}--\eqref{interface:condition:4}.
\end{lemma}

\section{An a priori estimate}
\label{sec:priori:estimate}

In this section, we shall propose an a priori estimate for possible solutions of \eqref{continuous:weak:formulation}.

\subsection{Some technique inequalities}

Firstly, throughout this paper we use $C$ to denote a generic positive constant independent of discretization parameters, which may take different values in different occasions. Then, based on general Sobolev inequalities, the trace theorem and the Sobolev embedding theorem \cite{Adams2003}, we have that for any bounded open set $\mathcal D\subset\mathbb{R}^N$ with Lipschitz continuous boundary $\partial\mathcal D$ and for all $v\in H^1(\mathcal D)$,

\begin{align}
\label{lq:domain}
\|v\|_{0,q,\mathcal D}&\leqslant C\|v\|_{1,\mathcal D},\quad 1\leqslant q\leqslant 6,\\
\label{l2:partial:domain}
\|v\|_{0,\partial\mathcal D}&\leqslant C\|v\|_{0,\mathcal D}^{1/2}\|v\|_{1,\mathcal D}^{1/2}\leqslant C\|v\|_{1,\mathcal D},\\
\label{l4:partial:domain}
\|v\|_{0,4,\partial\mathcal D}&\leqslant C\|v\|_{1,\mathcal D},
\end{align}
which can be also found in \cite{Chidyagwai20093806, Girault20092052}. We also need Poincar\'{e} inequality and Korn's inequality \cite{Cao20101} that for all $\bm{v}_s\in\bm{X}_s$ and for all $\psi\in X_d$,

\begin{align}
\label{v:Poincare}
\|\bm{v}_s\|_{0,\Omega_s}&\leqslant C|\bm{v}_s|_{1,\Omega_s},\\
\label{psi:Poincare}
\|\psi\|_{0,\Omega_d}&\leqslant C|\psi|_{1,\Omega_d},\\
\label{v:Korn}
|\bm{v}_s|_{1,\Omega_s}&\leqslant C\|\mathbb{D}(\bm{v}_s)\|_{0,\Omega_s}.
\end{align}
Moreover, for the trilinear form $b(\cdot;\cdot,\cdot)$, the following lemma holds with the help of \cite{Aubin1982, Bergh1976}, \eqref{lq:domain}, \eqref{v:Poincare} and \eqref{v:Korn}.

\begin{lemma}\label{lemma:trilinear:inequality}
For any functions $\bm{u}_s$, $\bm{v}_s$, $\bm{w}_s\in\bm{X}_s$, we have

\begin{align}
\label{trilinear:embedding}
\left|\left((\bm{u}_s\cdot\nabla)\bm{v}_s,\bm{w}_s\right)_{\Omega_s}\right|&\leqslant
C\|\bm{u}_s\|_{0,\Omega_s}^{1/2}\|\mathbb{D}(\bm{u}_s)\|_{0,\Omega_s}^{1/2}
\|\mathbb{D}(\bm{v}_s)\|_{0,\Omega_s}\|\mathbb{D}(\bm{w}_s)\|_{0,\Omega_s},\\
\label{trilinear:aux:embedding}
\left|\left((\nabla\cdot\bm{u}_s)\bm{v}_s,\bm{w}_s\right)_{\Omega_s}\right|&\leqslant
C\|\bm{v}_s\|_{0,\Omega_s}^{1/2}\|\mathbb{D}(\bm{v}_s)\|_{0,\Omega_s}^{1/2}
\|\mathbb{D}(\bm{u}_s)\|_{0,\Omega_s}\|\mathbb{D}(\bm{w}_s)\|_{0,\Omega_s}.
\end{align}
\end{lemma}

\begin{pf}
First, let us recall the standard interpolation inequality \cite{Aubin1982, Bergh1976}: for any bounded open set $\mathcal D\subset\mathbb{R}^N$ with Lipschitz continuous boundary $\partial\mathcal D$, $1\leqslant p<q<\infty$ and $\theta=N(1/p-1/q)\leqslant 1$,

\begin{equation}\label{interpolation:inequality}
  \|v\|_{0,q,\mathcal D}\leqslant C\|v\|_{0,p,\mathcal D}^{1-\theta}\|v\|_{1,p,\mathcal D}^\theta,
  \quad \forall\,v\in W^{1,p}(\mathcal D).
\end{equation}
In \eqref{interpolation:inequality}, we set $\mathcal D=\Omega_s$, and
\begin{itemize}
  \item if $N=2$, $p=2$, $q=4$, we then have

  \begin{equation}\label{2d:interpolation:inequality}
    \|\bm{v}_s\|_{0,4,\Omega_s}\leqslant C\|\bm{v}_s\|_{0,\Omega_s}^{1/2}\|\bm{v}_s\|_{1,\Omega_s}^{1/2},
    \quad \forall\,\bm{v}_s\in \bm{X}_s;
  \end{equation}
  \item if $N=3$, $p=2$, $q=3$, we then have

  \begin{equation}\label{3d:interpolation:inequality}
    \|\bm{v}_s\|_{0,3,\Omega_s}\leqslant C\|\bm{v}_s\|_{0,\Omega_s}^{1/2}\|\bm{v}_s\|_{1,\Omega_s}^{1/2},
    \quad \forall\,\bm{v}_s\in \bm{X}_s.
  \end{equation}
\end{itemize}
Hence, the results \eqref{trilinear:embedding} and \eqref{trilinear:aux:embedding} then follow from H\"{o}lder inequality, \eqref{lq:domain}, \eqref{v:Poincare}, \eqref{v:Korn}, \eqref{2d:interpolation:inequality} and \eqref{3d:interpolation:inequality} when $N=2$ or $3$.

\end{pf}

\subsection{The equivalent problem}

Let us introduce the following divergence-free space:

\[
\bm{V}_s:=\left\{\bm{v}_s\in\bm{X}_s:\,\nabla\cdot\bm{v}_s = 0\right\},
\]
and we consider the new product Hilbert space:

\[
\underline{\bm{W}}:=\bm{V}_s\times X_d\times X_d,
\]
which equipped with the same norm as $\underline{\bm{Y}}$.

Thanks to \cite{Girault20092052}, it provides a lower constant bound $\beta_0>0$ depending only on $\Omega$ to guarantee the following inf-sup condition:

\begin{equation}\label{con:infsup}
\inf_{q_s\in Q_s}\sup_{\underline{\bm{v}}=(\bm{v}_s,\psi_m,\psi_f)\in\underline{\bm{Y}}}
\frac{d(\bm{v}_s,q_s)}{\|q_s\|_{0,\Omega_s}\|\underline{\bm{v}}\|_{\underline{\bm{Y}}}}\geqslant\beta_0.
\end{equation}
Hence, by \eqref{con:infsup} and the same argument in \cite{Girault1986}, the Galerkin variational formulation \eqref{continuous:weak:formulation} is equivalent to the following problem: to find $\underline{\bm{u}}\in\underline{\bm{W}}$ such that

\begin{equation}\label{continuous:weak:formulation:2}
a(\underline{\bm{u}},\underline{\bm{v}})
+b(\bm{u}_s;\bm{u}_s,\bm{v}_s)
=\rho(\bm{f}_s,\bm{v}_s)_{\Omega_s}
+(f_d,\psi_f)_{\Omega_d},\quad
\forall\,\underline{\bm{v}}\in\underline{\bm{W}}.
\end{equation}

\subsection{Discretization and an auxiliary problem}

Let $\mathcal R_{i,h}$ be a quasi-uniform regular triangulation of the domain $\Omega_i$, $i=s,d$, respectively. If assuming that the two meshes $\mathcal R_{s,h}$ and $\mathcal R_{d,h}$ coincide along the interface $\Gamma$, then we can define $\mathcal R_h:=\mathcal R_{s,h}\cup\Gamma\cup\mathcal R_{d,h}$, which is also a quasi-uniform regular triangulation of $\Omega$. The diameter of element $K\in\mathcal R_h$ is denoted as $h_K$, and we set the mesh parameter $h=\max_{K\in\mathcal R_h}h_K$.

Then we denote by $\bm{X}_h\subset \bm{H}_0^1(\Omega)$ a finite element space defined on $\Omega$. The discrete space $\bm{X}_h$ can be naturally restricted to $\Omega_s$, so we set $\bm{X}_{s,h}=\bm{X}_h|_{\Omega_s}\subset\bm{X}_s$. Following the same technique, we establish other finite element spaces $Q_{s,h}\subset Q_s$ and $X_{d,h}\subset X_d$. We assume that $(\bm{X}_{s,h},Q_{s,h})$ is a stable finite element space pair. Further we define the following vector-valued Hilbert space on $\Omega_d$:

\[
\bm{X}_d:=\left\{\bm{v}\in\bm{H}^1(\Omega_d):\,\bm{v}=\bm{0}~\mathrm{on}~\Gamma_d\right\}.
\]
In addition, we need a finite element subspace $\bm{V}_h\subset\bm{X}_h$ defined on $\Omega$ given by

\[
\bm{V}_h:=\left\{\bm{v}_h\in\bm{X}_h:\,d(\bm{v}_h,q_{s,h})=0,\,\forall\,q_{s,h}\in Q_{s,h}\right\}.
\]
Similarly, we have

\[
\bm{V}_{s,h}=\bm{V}_h|_{\Omega_s}\subset\bm{X}_s,\quad \bm{X}_{d,h}=\bm{V}_h|_{\Omega_d}\subset\bm{X}_d,
\]
where $\bm{V}_{s,h}$ is a weakly divergence-free finite element space on $\Omega_s$, and any function $\bm{v}_{d,h}\in\bm{X}_{d,h}$ has an implicit restriction that $\int_\Gamma\bm{v}_{d,h}\cdot\bm{n}_d\mathrm{d}s = 0$ because of the continuity of $\bm{V}_h$ across the interface $\Gamma$.

According to \cite{Girault20092052}, the difficulty of obtaining an a priori estimate of the Navier-Stokes/Darcy fluid flow model comes from the energy unbalance due to the nonlinear convection term $(\bm{u}_s\cdot\nabla)\bm{u}_s$ in \eqref{stokes}. It motivates us to construct an auxiliary discrete Galerkin variational problem defined on $\Omega_d$ with compatible boundary conditions, so that the auxiliary problem can compensate the nonlinear convection term in the energy balance of the Navier-Stokes equations.

To fix ideas, we first define a lifting operator $\mathcal L:\,\bm{H}_{00}^{1/2}(\Gamma)\rightarrow\bm{X}_d$ as follows: for any $\bm{\eta}\in\bm{H}_{00}^{1/2}(\Gamma)$ with $\int_\Gamma\bm{\eta}\mathrm{d}s=0$ such that

\[
\mathcal L\bm{\eta}\in\bm{X}_d,\quad (\mathcal L\bm{\eta})|_\Gamma=\bm{\eta},\quad \nabla\cdot(\mathcal L\bm{\eta})|_{\Omega_d}=0.
\]
Then, we introduce the Scott-Zhang interpolator $\Pi_{s,h}:\,\bm{X}_s\rightarrow\bm{X}_{s,h}$ satisfying the following properties \cite{Scott1990483}:

\begin{align}
\label{SZ:error}
\|\bm{v}_s-\Pi_{s,h}\bm{v}_s\|_{0,\Omega_s}&\leqslant Ch\|\mathbb{D}(\bm{v}_s)\|_{0,\Omega_s},\quad\forall\,\bm{v}_s\in\bm{X}_s,\\
\label{SZ:h1:norm}
\|\Pi_{s,h}\bm{v}_s\|_{1,\Omega_s}&\leqslant
C\|\mathbb{D}(\bm{v}_s)\|_{0,\Omega_s},\quad\ \;\forall\,\bm{v}_s\in\bm{X}_s.
\end{align}
Now, with the mesh parameter $h$, the interpolator $\Pi_{s,h}$ and any given $\bm{u}_s\in\bm{V}_s$, we consider the following auxiliary discrete Galerkin variational problem: to find $\bm{u}_{d,h}\in\bm{X}_{d,h}$ with $\bm{u}_{d,h}|_\Gamma=(\Pi_{s,h}\bm{u}_s)|_\Gamma$ such that for all $\bm{v}_{d,h}\in\bm{X}_{d,h}$,

\begin{equation}\label{auxiliary:pde}
2\kappa\left(\mathbb{D}(\bm{u}_{d,h}),\mathbb{D}(\bm{v}_{d,h})\right)_{\Omega_d}
+\left((\bm{u}_d^0\cdot\nabla)\bm{u}_{d,h},\bm{v}_{d,h}\right)_{\Omega_d}-
\kappa\left\langle\frac{\partial\bm{u}_{d,h}}{\partial\bm{n}_d},\bm{v}_{d,h}\right\rangle_\Gamma=0,
\end{equation}
where $\bm{u}_d^0=\mathcal L(\bm{u}_s|_\Gamma)\in\bm{X}_d$, and $\kappa>0$ is a certain positive constant specified later. Furthermore, the discrete problem \eqref{auxiliary:pde} is uniquely solvable in $\bm{X}_{d,h}$ from Lemma 3.1 in \cite{Hou201947}.

\begin{remark}
The discrete problem \eqref{auxiliary:pde} is the conforming Galerkin approximation of the following convection-diffusion equation defined on $\Omega_d$: for any given $\bm{u}_s\in\bm{V}_s$, to find $\bm{u}_d\in\bm{X}_d$ with $\bm{u}_d|_\Gamma=\bm{u}_s|_\Gamma$ such that

\begin{equation}\label{auxiliary:pde:con}
  \begin{cases}
    -2\kappa\nabla\cdot\mathbb{D}(\bm{u}_d)+(\bm{u}_d^0\cdot\nabla)\bm{u}_d=\bm{0}
    & \mathrm{in}~\Omega_d,\\
    \bm{u}_d=\bm{0} & \mathrm{on}~\Gamma_d.
  \end{cases}
\end{equation}
Clearly, for any constant $\kappa>0$ and $\bm{u}_s\in\bm{V}_s$, the solution of \eqref{auxiliary:pde:con} is well-posed \cite{Roos2008}.
\end{remark}

\subsection{An a priori estimate of weak solutions}\label{subsec:a:priori}

We realize that \eqref{continuous:weak:formulation:2} and \eqref{auxiliary:pde} can form a larger coupled system because of the convection term $\bm{u}_d^0$ and the constraint $\bm{u}_{d,h}|_\Gamma=(\Pi_{s,h}\bm{u}_s)|_\Gamma$ for the unknown variable $\bm{u}_{d,h}$ over the interface $\Gamma$, and we also stress that the new coupled system has no energy exchange via the interface $\Gamma$. It implies that \eqref{auxiliary:pde} is subjected to \eqref{continuous:weak:formulation:2} but any possible solution of \eqref{continuous:weak:formulation:2} has nothing to do with the auxiliary problem \eqref{auxiliary:pde}.

\begin{theorem}\label{theorem:priori:estimate}
Assume that the data in the auxiliary discrete problem \eqref{auxiliary:pde} satisfy $h$ small enough and $0<\kappa\leqslant C\rho\nu h$. If problem \eqref{continuous:weak:formulation:2} exists a possible solution $(\bm{u}_s,\phi_m,\phi_f)\in\underline{\bm{W}}$, we then have the following a priori estimate that

\begin{equation}\label{a:priori:estimate}
  \rho\nu\|\mathbb{D}(\bm{u}_s)\|_{0,\Omega_s}^2+\frac{2\sigma k_m}{\mu}\|\nabla\phi_m\|_{0,\Omega_d}^2
  +\frac{\sigma k_f}{\mu}\|\nabla\phi_f\|_{0,\Omega_d}^2\leqslant\mathcal C^2,
\end{equation}
where

\[
\mathcal C^2=\rho\nu^{-1}\|\bm{f}_s\|_{-1,\Omega_s}^2
  +\mu \sigma^{-1} k_f^{-1}\|f_d\|_{-1,\Omega_d}^2.
\]
Here we stress that there are no assumptions for the data and physical parameters of the model problem.
\end{theorem}

\begin{pf}

We denote by $\underline{\bm{u}}=(\bm{u}_s,\phi_m,\phi_f)\in\underline{\bm{W}}$ a possible solution of problem \eqref{continuous:weak:formulation:2}, and denote by $\bm{u}_{d,h}\in\bm{X}_{d,h}$ the solution of problem \eqref{auxiliary:pde}. Then, we assume that there is a positive finite constant $M_s<\infty$ such that $\|\mathbb{D}(\bm{u}_s)\|_{0,\Omega_s}\leqslant M_s$. Taking $\underline{\bm{v}}=\underline{\bm{u}}$ in \eqref{continuous:weak:formulation:2}, and noting that the terms $a_\Gamma(\underline{\bm{u}},\underline{\bm{u}})$ and $\frac{\sigma k_m}{\mu}\left(\phi_m-\phi_f,\phi_m\right)_{\Omega_d}+\frac{\sigma k_m}{\mu}\left(\phi_f-\phi_m,\phi_f\right)_{\Omega_d}$ are non-negative, we have

\begin{equation}\label{priori:estimate:temp1}
  2\rho\nu\|\mathbb{D}(\bm{u}_s)\|_{0,\Omega_s}^2+\frac{\sigma k_m}{\mu}\|\nabla\phi_m\|_{0,\Omega_d}^2
  +\frac{\sigma k_f}{\mu}\|\nabla\phi_f\|_{0,\Omega_d}^2+\frac{1}{2}\left\langle\bm{u}_s\cdot\bm{n}_s,|\bm{u}_s|^2\right\rangle_\Gamma\leqslant \rho(\bm{f}_s,\bm{u}_s)_{\Omega_s}+(f_d,\phi_f)_{\Omega_d}.
\end{equation}
In addition, we note that $\nabla\cdot\bm{u}_d^0=0$ in $\Omega_d$, and the $\bm{u}_d^0|_\Gamma=\bm{u}_s|_\Gamma$ by the definition of $\mathcal L$. Hence, the identity \eqref{trilinear:form} can also be applied to the second term in the left hand of \eqref{auxiliary:pde}, that is for all $\bm{v}_{d,h}\in\bm{X}_{d,h}$,

\begin{equation}\label{priori:estimate:temp2}
\begin{split}
&\left((\bm{u}_d^0\cdot\nabla)\bm{u}_{d,h},\bm{v}_{d,h}\right)_{\Omega_d}
=\left((\bm{u}_d^0\cdot\nabla)\bm{u}_{d,h},\bm{v}_{d,h}\right)_{\Omega_d}
+\frac{1}{2}\left((\nabla\cdot\bm{u}_d^0)\bm{u}_{d,h},\bm{v}_{d,h}\right)_{\Omega_d}\\
&=-\frac{1}{2}\left\langle\bm{u}_s\cdot\bm{n}_s,|\Pi_{s,h}\bm{u}_s|^2\right\rangle_\Gamma
+\frac{1}{2}\left((\bm{u}_d^0\cdot\nabla)\bm{u}_{d,h},\bm{v}_{d,h}\right)_{\Omega_d}
-\frac{1}{2}\left((\bm{u}_d^0\cdot\nabla)\bm{v}_{d,h},\bm{u}_{d,h}\right)_{\Omega_d}.
\end{split}
\end{equation}
So, taking $\bm{v}_{d,h}=\bm{u}_{d,h}$ in \eqref{auxiliary:pde} and using \eqref{priori:estimate:temp2}, we obtain

\begin{equation}\label{priori:estimate:temp3}
  2\kappa\|\mathbb{D}(\bm{u}_{d,h})\|_{0,\Omega_d}^2
  -\frac{1}{2}\left\langle\bm{u}_s\cdot\bm{n}_s,|\Pi_{s,h}\bm{u}_s|^2\right\rangle_\Gamma
  =\kappa\left\langle\frac{\partial\bm{u}_{d,h}}{\partial\bm{n}_d},\Pi_{s,h}\bm{u}_s\right\rangle_\Gamma.
\end{equation}

It follows from \eqref{l2:partial:domain}, \eqref{l4:partial:domain}, \eqref{v:Poincare}, \eqref{v:Korn}, \eqref{SZ:error}, \eqref{SZ:h1:norm}, H\"{o}lder inequality and the triangle inequality that

\begin{equation}\label{priori:estimate:temp4}
\begin{split}
&\frac{1}{2}\left\langle\bm{u}_s\cdot\bm{n}_s,|\bm{u}_s|^2\right\rangle_\Gamma
-\frac{1}{2}\left\langle\bm{u}_s\cdot\bm{n}_s,|\Pi_{s,h}\bm{u}_s|^2\right\rangle_\Gamma\\
\leqslant&
\frac{1}{2}\left\|\bm{u}_s\right\|_{0,4,\Gamma}\|\bm{u}_s+\Pi_{s,h}\bm{u}_s\|_{0,4,\Gamma}\left\|\bm{u}_s-\Pi_{s,h}\bm{u}_s\right\|_{0,\Gamma}\\
\leqslant& Ch^{1/2}M_s\|\mathbb{D}(\bm{u}_s)\|_{0,\Omega_s}^2.
\end{split}
\end{equation}
As follows, we now define the dual norms of $\bm{X}_s$ and $X_d$, which are denoted as $\|\cdot\|_{-1,\Omega_s}$ and $\|\cdot\|_{-1,\Omega_d}$, respectively.

\[
\|\bm{f}_s\|_{-1,\Omega_s}:=\inf_{\bm{v}_s\in\bm{X}_s}\frac{(\bm{f}_s,\bm{v}_s)_{\Omega_s}}{\|\mathbb{D}(\bm{v}_s)\|_{0,\Omega_s}},
\quad \|f_d\|_{-1,\Omega_d}:=\inf_{\psi\in X_d}\frac{(f_d,\psi)_{\Omega_d}}{\|\nabla\psi\|_{0,\Omega_d}}.
\]
Hence, by using H\"{o}lder inequality and Young's inequality, we have

\begin{equation}\label{priori:estimate:temp5}
  \rho(\bm{f}_s,\bm{u}_s)_{\Omega_s}+(f_d,\phi_f)_{\Omega_d}
  \leqslant \frac{\rho\nu}{2}\|\mathbb{D}(\bm{u}_s)\|_{0,\Omega_s}^2
  +\frac{\sigma k_f}{2\mu}\|\nabla\phi_f\|_{0,\Omega_d}^2+\frac{\rho}{2\nu}\|\bm{f}_s\|_{-1,\Omega_s}^2
  +\frac{\mu}{2\sigma k_f}\|f_d\|_{-1,\Omega_d}^2.
\end{equation}
In addition, based on the imbedding result \eqref{interface:space:imbedding}, the standard trace theorem \cite{Adams2003}, H\"{o}lder inequality, Korn's inequality, Young's inequality, \eqref{SZ:h1:norm} and the following inverse inequality \cite{Riviere2008}: for any polynomial $\bm{v}$ on $K$,

\[
\|(\nabla\bm{v})\bm{n}\|_{0,e}\leqslant Ch^{-1/2}|\bm{v}|_{1,K},\quad\forall\,e\subset\partial K,\ \forall\,K\in\mathcal R_h,
\]
we have

\begin{equation}\label{priori:estimate:temp6}
\begin{split}
  &\kappa\left\langle\frac{\partial\bm{u}_{d,h}}{\partial\bm{n}_d},\Pi_{s,h}\bm{u}_s\right\rangle_\Gamma
  \leqslant \kappa\left\|(\nabla\bm{u}_{d,h})\bm{n}_d\right\|_{\left(\bm{H}_{00}^{1/2}(\Gamma)\right)'}
  \left\|\Pi_{s,h}\bm{u}_s\right\|_{\bm{H}_{00}^{1/2}(\Gamma)}\\
  &\quad\leqslant C\kappa
  \left\|(\nabla\bm{u}_{d,h})\bm{n}_d\right\|_{0,\Gamma}\left\|\Pi_{s,h}\bm{u}_s\right\|_{1,\Omega_s}
  \leqslant C\kappa
  \left(\sum_{e\in\Gamma}\|(\nabla\bm{u}_{d,h})\bm{n}\|_{0,e}^2\right)^{1/2}
  \left\|\mathbb{D}(\bm{u}_s)\right\|_{0,\Omega_s}\\
  &\quad\leqslant C\kappa h^{-1/2}\left\|\mathbb{D}(\bm{u}_{d,h})\right\|_{0,\Omega_d}
  \left\|\mathbb{D}(\bm{u}_s)\right\|_{0,\Omega_s}\leqslant
  \frac{C\kappa^2}{\rho\nu h}\left\|\mathbb{D}(\bm{u}_{d,h})\right\|_{0,\Omega_d}^2
  +\frac{\rho\nu}{2}\left\|\mathbb{D}(\bm{u}_s)\right\|_{0,\Omega_s}^2.
\end{split}
\end{equation}
Finally, with the assumptions that $h$ is small enough such that $Ch^{1/2}M_s<\rho\nu/2$, and $\kappa$ satisfies $0<\kappa\leqslant C\rho\nu h$, gathering \eqref{priori:estimate:temp1}, \eqref{priori:estimate:temp3}, \eqref{priori:estimate:temp4}, \eqref{priori:estimate:temp5} and \eqref{priori:estimate:temp6} yields \eqref{a:priori:estimate}, where

\[
\mathcal C^2=\rho\nu^{-1}\|\bm{f}_s\|_{-1,\Omega_s}^2
  +\mu \sigma^{-1} k_f^{-1}\|f_d\|_{-1,\Omega_d}^2.
\]
\end{pf}

\section{Existence and global uniqueness of the solution}

In this section, we shall use the technique of the Galerkin method to verify that problem \eqref{continuous:weak:formulation:2} has at least one weak solution, and then we can prove the global uniqueness of the weak solution due to the a priori estimate \eqref{a:priori:estimate} obtained in Section \ref{sec:priori:estimate}.

\subsection{The solvability of the conforming Galerkin approximation problem}

We denote by $\underline{\bm{W}}_h$ the product finite element space $\bm{V}_{s,h}\times X_{d,h}\times X_{d,h}$. Then, we consider the following conforming Galerkin approximation problem of \eqref{continuous:weak:formulation:2}: to find $\underline{\bm{u}}_h=(\bm{u}_{s,h},\phi_{m,h},\phi_{f,h})\in\underline{\bm{W}}_h$ such that $\forall\,\underline{\bm{v}}_h=(\bm{v}_{s,h},\psi_{m,h},\psi_{f,h})\in\underline{\bm{W}}_h$,

\begin{equation}\label{discrete:weak:formulation}
a(\underline{\bm{u}}_h,\underline{\bm{v}}_h)
+b(\bm{u}_{s,h};\bm{u}_{s,h},\bm{v}_{s,h})
=\rho(\bm{f}_s,\bm{v}_{s,h})_{\Omega_s}
+(f_d,\psi_{f,h})_{\Omega_d}.
\end{equation}
To show the solvability of \eqref{discrete:weak:formulation}, stemming from Theorem \ref{theorem:priori:estimate}, we shall consider the following constructed coupled discrete system: to find $(\underline{\bm{u}}_h,\bm{u}_{d,h})\in\underline{\bm{W}}_h\times \bm{X}_{d,h}$ such that

\begin{align}
\label{discrete:au:1}
&a(\underline{\bm{u}}_h,\underline{\bm{v}}_h)
+b(\bm{u}_{s,h};\bm{u}_{s,h},\bm{v}_{s,h})
=\rho(\bm{f}_s,\bm{v}_{s,h})_{\Omega_s}
+(f_d,\psi_{f,h})_{\Omega_d}
& \forall\,\underline{\bm{v}}_h\in\underline{\bm{W}}_h,\\
\label{discrete:au:2}
&2\xi\left(\mathbb{D}(\bm{u}_{d,h}),\mathbb{D}(\bm{v}_{d,h})\right)_{\Omega_d}
+\left((\bm{\beta}\cdot\nabla)\bm{u}_{d,h},\bm{v}_{d,h}\right)_{\Omega_d}
-\xi\left\langle\frac{\partial\bm{u}_{d,h}}{\partial\bm{n}_d},\bm{v}_{d,h}\right\rangle_\Gamma
=0&\forall\,\bm{v}_{d,h}\in\bm{X}_{d,h},\\
\label{discrete:au:3}
&\bm{u}_{d,h}|_\Gamma=\bm{u}_{s,h}|_\Gamma,&
\end{align}
where $\xi>0$ is a positive constant specified later, $\bm{\beta}:=\mathcal L(\bm{u}_{s,h}|_\Gamma)$, and $\nabla\cdot\bm{\beta}=0$ in $\Omega_d$ by the definition of the lifting operator $\mathcal L$. Furthermore, it follows the head statements of Section \ref{subsec:a:priori} that if $(\underline{\bm{u}}_h,\bm{u}_{d,h})\in\underline{\bm{W}}_h\times \bm{X}_{d,h}$ is a solution of problem \eqref{discrete:au:1}--\eqref{discrete:au:3}, then $\underline{\bm{u}}_h\in\underline{\bm{W}}_h$ will solve \eqref{discrete:weak:formulation}.

\begin{lemma}\label{lemma:discrete:priori:estimate}
If problem \eqref{discrete:weak:formulation} exists a possible solution $\underline{\bm{u}}_h=(\bm{u}_{s,h},\phi_{m,h},\phi_{f,h})\in\underline{\bm{W}}_h$, we have the following a priori estimate that

\begin{equation}\label{discrete:a:priori:estimate}
  \rho\nu\|\mathbb{D}(\bm{u}_{s,h})\|_{0,\Omega_s}^2+\frac{2\sigma k_m}{\mu}\|\nabla\phi_{m,h}\|_{0,\Omega_d}^2
  +\frac{\sigma k_f}{\mu}\|\nabla\phi_{f,h}\|_{0,\Omega_d}^2\leqslant\mathcal C^2,
\end{equation}
where $\mathcal C^2$ is defined in Theorem \ref{theorem:priori:estimate}.
\end{lemma}

\begin{pf}
The proof is quite close to the proof of Theorem \ref{theorem:priori:estimate}. To avoid repeating, we just present the differences. Since $\bm{\beta}=\mathcal L(\bm{u}_{s,h}|_\Gamma)$, we have

\[
\left((\bm{\beta}\cdot\nabla)\bm{u}_{d,h},\bm{u}_{d,h}\right)_{\Omega_d}
=-\frac{1}{2}\left\langle\bm{u}_{s,h}\cdot\bm{n}_s,\left|\bm{u}_{s,h}\right|^2\right\rangle_\Gamma,
\]
and thus the term \eqref{priori:estimate:temp4} vanishes here. As a result, some subtle differences occur to the following estimate:

\[
\xi\left\langle\frac{\partial\bm{u}_{d,h}}{\partial\bm{n}_d},\bm{u}_{s,h}\right\rangle_\Gamma
\leqslant \frac{C\xi^2}{2\rho\nu h}\|\mathbb{D}(\bm{u}_{d,h})\|_{0,\Omega_d}^2
+\rho\nu\|\mathbb{D}(\bm{u}_{s,h})\|_{0,\Omega_s}^2.
\]
Thus, for any given mesh parameter $h$ , the result \eqref{discrete:a:priori:estimate} follows the assumption that $0<\xi\leqslant C\rho\nu h$.

\end{pf}

Now we start to verify the solvability of problem \eqref{discrete:au:1}--\eqref{discrete:au:3}. For any $\widehat{\bm{v}}_h\in\bm{V}_h$, following the similar technique proposed in \cite{Hou201947}, we denote $\bm{v}_h^s=\widehat{\bm{v}}_h|_{\Omega_s}\in\bm{V}_{s,h}$,  $\bm{v}_h^d=\widehat{\bm{v}}_h|_{\Omega_d}\in\bm{X}_{d,h}$, and

\[
\widehat{\widetilde{\bm{v}}}_h:=\left\{
\begin{array}{cl}
\bm{v}_h^s&\mathrm{in}~\Omega_s,\\
\bm{v}_d^0&\mathrm{in}~\Omega_d,
\end{array}
\right.
\]
where $\bm{v}_d^0=\mathcal L(\bm{v}_h^s|_\Gamma)\in\bm{X}_d$ with $\nabla\cdot\bm{v}_d^0=0$ in the domain $\Omega_d$. Then, we denote by $\underline{\bm{Z}}_h$ the product space $\bm{V}_h\times X_{d,h}\times X_{d,h}$, and define a mapping $\mathcal F_h:\underline{\bm{Z}}_h\rightarrow\underline{\bm{Z}}_h$ as: for each  $\widehat{\underline{\bm{v}}}_h=(\widehat{\bm{v}}_h,\psi_{m,h},\psi_{f,h})\in \underline{\bm{Z}}_h$ such that $\forall\,\widehat{\underline{\bm{w}}}_h=(\widehat{\bm{w}}_h,\varphi_{m,h},\varphi_{f,h})\in \underline{\bm{Z}}_h$,

\begin{equation}\label{discrete:mapping}
\begin{split}
&\left(\mathcal F_h\widehat{\underline{\bm{v}}}_h,\widehat{\underline{\bm{w}}}_h\right)_{\underline{\bm{Z}}_h}
:=2\rho\nu\left(\mathbb{D}(\bm{v}_h^s),\mathbb{D}(\bm{w}_h^s)\right)_{\Omega_s}
+\frac{k_m}{\mu}\left(\nabla\psi_{m,h},\nabla\varphi_{m,h}\right)_{\Omega_d}+
\frac{k_f}{\mu}\left(\nabla\psi_{f,h},\nabla\varphi_{f,h}\right)_{\Omega_d}\\
&\quad+\frac{\sigma k_m}{\mu}\left(\psi_{m,h}-\psi_{f,h},\varphi_{m,h}\right)_{\Omega_d}
+\frac{\sigma k_m}{\mu}\left(\psi_{f,h}-\psi_{m,h},\varphi_{f,h}\right)_{\Omega_d}
+2\xi\left(\mathbb{D}(\bm{v}_h^d),\mathbb{D}(\bm{w}_h^d)\right)_{\Omega_d}\\
&\quad
+\left\langle\psi_{f,h},\bm{w}_h^s\cdot\bm{n}_s\right\rangle_\Gamma
-\left\langle\varphi_{f,h},\bm{v}_h^s\cdot\bm{n}_s\right\rangle_\Gamma
+\sum_{i=1}^{N-1}\frac{\alpha\rho\nu}{\sqrt{k_f}}
\left\langle\bm{v}_h^s\cdot\bm{\tau}_i,\bm{w}_h^s\cdot\bm{\tau}_i\right\rangle_\Gamma
-\xi\left\langle\frac{\partial\bm{v}_h^d}{\partial\bm{n}_d},\bm{w}_h^d\right\rangle_\Gamma\\
&\quad+\left((\widehat{\widetilde{\bm{v}}}_h\cdot\nabla)\widehat{\bm{v}}_h,\widehat{\bm{w}}_h\right)_{\Omega}
+\frac{1}{2}\left((\nabla\cdot\widehat{\widetilde{\bm{v}}}_h)\widehat{\bm{v}}_h,\widehat{\bm{w}}_h\right)_{\Omega}
-\rho(\bm{f}_s,\bm{w}_h^s)_{\Omega_s}-(f_d,\varphi_{f,h})_{\Omega_d}.
\end{split}
\end{equation}
Clearly, $\mathcal F_h$ define a mapping from $\underline{\bm{Z}}_h$ into itself, and a zero of $\mathcal F_h$ is a solution of the coupled system \eqref{discrete:au:1}--\eqref{discrete:au:3}. Further we introduce the Brouwer's fixed-point theorem:

\begin{lemma}\label{lemma:Brouwer:fix}\cite{Cioranescu2016}
Let $U$ be a nonempty, convex, and compact subset of a normed vector space and let $F$ be a continuous mapping from $U$ into $U$. Then $F$ has at least one fixed point.
\end{lemma}
Based on Lemma \ref{lemma:discrete:priori:estimate}, we define $\underline{\bm{U}}_h$ a subset of $\underline{\bm{Z}}_h$ as:

\[
\underline{\bm{U}}_h:=\left\{\widehat{\underline{\bm{v}}}_h\in\underline{\bm{Z}}_h:\,
\rho\nu\|\mathbb{D}(\bm{v}_h^s)\|_{0,\Omega_s}^2+\frac{2\sigma k_m}{\mu}\|\nabla\psi_{m,h}\|_{0,\Omega_d}^2+\frac{\sigma k_f}{\mu}\|\nabla\psi_{f,h}\|_{0,\Omega_d}^2\leqslant\mathcal C^2\right\},
\]
where $\mathcal C^2$ is defined in Theorem \ref{theorem:priori:estimate}. Then, taking $\widehat{\underline{\bm{w}}}_h=\widehat{\underline{\bm{v}}}_h\in\underline{\bm{U}}_h$ in \eqref{discrete:mapping}, and following the steps of proving Theorem \ref{theorem:priori:estimate} and Lemma \ref{lemma:discrete:priori:estimate}, we obtain that

\[
\left(\mathcal F_h\widehat{\underline{\bm{v}}}_h,
\widehat{\underline{\bm{v}}}_h\right)_{\underline{\bm{Z}}_h}\geqslant 0,
\quad\forall\,\widehat{\underline{\bm{v}}}_h\in\underline{\bm{U}}_h.
\]
Hence, it follows Lemma \ref{lemma:Brouwer:fix} that there is at least one zero of $\mathcal F_h$ in the ball $\underline{\bm{U}}_h$ centered at the origin.

Gathering all the above results, we conclude the following theorem:

\begin{theorem}\label{theorem:discrete:existence}
For any given mesh parameter $h>0$, the conforming Galerkin approximation problem \eqref{discrete:weak:formulation} has at least one solution $(\bm{u}_{s,h},\phi_{m,h},\phi_{f,h})\in\bm{V}_{s,h}\times X_{d,h}\times X_{d,h}$ with the following estimate:

\begin{equation}\label{discrete:a:priori:estimate:2}
  \rho\nu\|\mathbb{D}(\bm{u}_{s,h})\|_{0,\Omega_s}^2+\frac{2\sigma k_m}{\mu}\|\nabla\phi_{m,h}\|_{0,\Omega_d}^2
  +\frac{\sigma k_f}{\mu}\|\nabla\phi_{f,h}\|_{0,\Omega_d}^2\leqslant\mathcal C^2,
\end{equation}
where $\mathcal C^2$ is defined in Theorem \ref{theorem:priori:estimate}.
\end{theorem}

\subsection{Existence and global uniqueness}

It stems from Theorem \ref{theorem:discrete:existence} and the conforming property $\bm{V}_{s,h}\times\bm{X}_{d,h}\times\bm{X}_{d,h}\subset\bm{V}_s\times\bm{X}_d\times\bm{X}_d$ that there exists a 3-tuple function $(\bm{u}_s,\phi_m,\phi_f)$ in the Hilbert space $\underline{\bm{W}}=\bm{V}_s\times\bm{X}_d\times\bm{X}_d$, and a uniformly bounded subsequence $\left\{(\bm{u}_{s,h},\phi_{m,h},\phi_{f,h})\right\}_{h>0}$ such that

\begin{equation}\label{limit:h1}
\lim_{h\rightarrow 0}(\bm{u}_{s,h},\phi_{m,h},\phi_{f,h})=(\bm{u}_{s},\phi_m,\phi_f)\quad\mathrm{weakly~in}~\underline{\bm{W}}.
\end{equation}
%\begin{align}
%\label{limit:h1}
%\lim_{h\rightarrow 0}(\bm{u}_{s,h},\phi_{m,h},\phi_{f,h})&=(\bm{u}_{s},\phi_m,\phi_f)\quad\mathrm{weakly~in}~\underline{\bm{W}},\\
%\label{limit:h12:u}
%\lim_{h\rightarrow 0}\bm{u}_{s,h}|_{\partial\Omega_s}&=\bm{u}_s|_{\partial\Omega_s}\qquad\quad\ \mathrm{weakly~in}~\bm{H}^{1/2}(\partial\Omega_s),\\
%\label{limit:h12:phi}
%\lim_{h\rightarrow 0}\phi_{i,h}|_{\partial\Omega_d}&=\phi_i|_{\partial\Omega_d}\qquad\quad\ \, \mathrm{weakly~in}~H^{1/2}(\partial\Omega_d),\quad i=m,f.
%\end{align}
Furthermore, the Sobolev imbedding implies that the above convergence results are strong in $L^q(\Omega)$ for any $1\leqslant q<6$ whenever $N=2$ or $3$. In particular, by extracting another subsequence, still denoted by $h$, we obtain

\begin{equation}\label{limit:l2}
\lim_{h\rightarrow 0}(\bm{u}_{s,h},\phi_{m,h},\phi_{f,h})=(\bm{u}_{s},\phi_m,\phi_f)\quad\mathrm{strongly~in}~\bm{L}^2(\Omega_s)\times L^2(\Omega_d)\times L^2(\Omega_d).
\end{equation}

\begin{theorem}\label{theorem:unique:1}
If the data satisfies that

\begin{equation}\label{assum:unique}
\mathcal N\left(\rho^{-1}\nu^{-2}\|\bm{f}_s\|_{-1,\Omega_s}+\rho^{-3/2}\nu^{-3/2}\mu^{1/2}\sigma^{-1/2}k_f^{-1/2}\|f_d\|_{-1,\Omega_d}\right)<1,
\end{equation}
the problem \eqref{continuous:weak:formulation:2} then admits a unique solution $(\bm{u}_s,\phi_m,\phi_f)$ in $\bm{V}_s\times X_d\times X_d$ such that

\[
  \rho\nu\|\mathbb{D}(\bm{u}_{s})\|_{0,\Omega_s}^2+\frac{2\sigma k_m}{\mu}\|\nabla\phi_{m}\|_{0,\Omega_d}^2
  +\frac{\sigma k_f}{\mu}\|\nabla\phi_{f}\|_{0,\Omega_d}^2\leqslant\mathcal C^2,
\]
where $\mathcal C^2$ is defined in Theorem \ref{theorem:priori:estimate}.
\end{theorem}

\begin{pf}
For $(\bm{u}_{s},\phi_m,\phi_f)$ denoted as $\underline{\bm{u}}$ and the subsequences $\{(\bm{u}_{s,h},\phi_{m,h},\phi_{f,h})\}_{h>0}$ denoted as $\{\underline{\bm{u}}_h\}_{h>0}$ defined in \eqref{limit:h1} and \eqref{limit:l2}, we can easily obtain

\begin{equation}\label{limit:h:temp1}
\lim_{h\rightarrow 0}\{a_s(\underline{\bm{u}}_h,\underline{\bm{v}})
+a_d(\underline{\bm{u}}_h,\underline{\bm{v}})\}=a_s(\underline{\bm{u}},\underline{\bm{v}})
+a_d(\underline{\bm{u}},\underline{\bm{v}}),\quad\forall\,\underline{\bm{v}}=(\bm{v}_s,\psi_m,\psi_f)\in\underline{\bm{W}},
\end{equation}
because of \eqref{limit:h1}. Then, for the trace bilinear term $a_\Gamma(\underline{\bm{u}}_h,\underline{\bm{v}})$, it follows from H\"{o}lder inequality, \eqref{l2:partial:domain}, and \eqref{v:Poincare}--\eqref{v:Korn} that

\begin{align*}
&a_\Gamma(\underline{\bm{u}}_h,\underline{\bm{v}})\\
&=\left\langle\phi_{f,h}-\phi_f,\bm{v}_s\cdot\bm{n}_s\right\rangle_\Gamma
-\left\langle\psi_f,(\bm{u}_{s,h}-\bm{u}_s)\cdot\bm{n}_s\right\rangle_\Gamma
+\sum_{i=1}^{N-1}\frac{\alpha\rho\nu}{\sqrt{k_f}}
\left\langle(\bm{u}_{s,h}-\bm{u}_s)\cdot\bm{\tau}_i,\bm{v}_s\cdot\bm{\tau}_i\right\rangle_\Gamma\\
&\quad+\left\langle\phi_f,\bm{v}_s\cdot\bm{n}_s\right\rangle_\Gamma
-\left\langle\psi_f,\bm{u}_s\cdot\bm{n}_s\right\rangle_\Gamma
+\sum_{i=1}^{N-1}\frac{\alpha\rho\nu}{\sqrt{k_f}}
\left\langle\bm{u}_s\cdot\bm{\tau}_i,\bm{v}_s\cdot\bm{\tau}_i\right\rangle_\Gamma\\
&\leqslant C\|\phi_{f,h}-\phi_f\|_{0,\Omega_d}^{1/2}\|\nabla(\phi_{f,h}-\phi_f)\|_{0,\Omega_d}^{1/2}\|\mathbb{D}(\bm{v}_s)\|_{0,\Omega_s}
+C\|\bm{u}_{s,h}-\bm{u}_s\|_{0,\Omega_s}^{1/2}\|\mathbb{D}(\bm{u}_{s,h}-\bm{u}_s)\|_{0,\Omega_s}^{1/2}\|\nabla\psi_f\|_{0,\Omega_d}\\
&\quad+C\|\bm{u}_{s,h}-\bm{u}_s\|_{0,\Omega_s}^{1/2}\|\mathbb{D}(\bm{u}_{s,h}-\bm{u}_s)\|_{0,\Omega_s}^{1/2}\|\mathbb{D}(\bm{v}_s)\|_{0,\Omega_s}
+a_\Gamma(\underline{\bm{u}},\underline{\bm{v}}).
\end{align*}
Since the uniform boundedness of $(\bm{u}_{s,h},\phi_{m,h},\phi_{f,h})$ shown in \eqref{discrete:a:priori:estimate:2} and the strong $L^2$-convergence result \eqref{limit:l2}, we derive that

\begin{equation}\label{limit:h:temp2}
\lim_{h\rightarrow 0}a_\Gamma(\underline{\bm{u}}_h,\underline{\bm{v}})=a_\Gamma(\underline{\bm{u}},\underline{\bm{v}}),
\quad\forall\,\underline{\bm{v}}\in\underline{\bm{W}}.
\end{equation}
In addition, for the limit of the trilinear form $b(\bm{u}_{s,h};\bm{u}_{s,h},\bm{v}_s)$, gathering \eqref{trilinear:form} and Lemma \ref{lemma:trilinear:inequality} yields

\begin{align*}
&\left|b(\bm{u}_{s,h};\bm{u}_{s,h},\bm{v}_s)-b(\bm{u}_s;\bm{u}_s,\bm{v}_s)\right|\\
&\leqslant
\left|\left((\bm{u}_{s,h}-\bm{u}_s)\cdot\nabla)\bm{u}_{s,h},\bm{v}_s\right)_{\Omega_s}\right|
+\frac{1}{2}\left|\left(\nabla\cdot\bm{u}_{s,h},(\bm{u}_{s,h}-\bm{u}_s)\cdot\bm{v}_s\right)_{\Omega_s}\right|\\
&\quad
+\left|\left((\bm{u}_s\cdot\nabla)(\bm{u}_{s,h}-\bm{u}_s),\bm{v}_s\right)_{\Omega_s}\right|
+\frac{1}{2}\left|\left(\nabla\cdot(\bm{u}_{s,h}-\bm{u}_s),\bm{u}_s\cdot\bm{v}_s\right)_{\Omega_s}\right|\\
&\leqslant
\underbrace{C\|\bm{u}_{s,h}-\bm{u}_s\|_{0,\Omega_s}^{1/2}\|\mathbb{D}(\bm{u}_{s,h}-\bm{u}_s)\|_{0,\Omega_s}^{1/2}
\|\mathbb{D}(\bm{u}_{s,h})\|_{0,\Omega_s}\|\mathbb{D}(\bm{v}_s)\|_{0,\Omega_s}}_{\uppercase\expandafter{\romannumeral1}_h}\\
&\quad
+\underbrace{\left|\left((\bm{u}_s\cdot\nabla)(\bm{u}_{s,h}-\bm{u}_s),\bm{v}_s\right)_{\Omega_s}\right|
+\frac{1}{2}\left|\left(\nabla\cdot(\bm{u}_{s,h}-\bm{u}_s),\bm{u}_s\cdot\bm{v}_s\right)_{\Omega_s}\right|}_{\uppercase\expandafter{\romannumeral2}_h}.
\end{align*}
We can easily obtain that $\lim_{h\rightarrow 0}\uppercase\expandafter{\romannumeral1}_h=0$ by \eqref{discrete:a:priori:estimate:2} and \eqref{limit:l2}, and $\lim_{h\rightarrow 0}\uppercase\expandafter{\romannumeral2}_h=0$ by \eqref{limit:h1}. Therefore the following limit holds:

\begin{equation}\label{limit:h:temp3}
\lim_{h\rightarrow 0}b(\bm{u}_{s,h};\bm{u}_{s,h},\bm{v}_s)=b(\bm{u}_s;\bm{u}_s,\bm{v}_s),
\quad\forall\,\bm{v}_s\in\bm{V}_s.
\end{equation}
It follows from \eqref{limit:h:temp1}--\eqref{limit:h:temp3} that

\[
a(\underline{\bm{u}},\underline{\bm{v}})
+b(\bm{u}_s;\bm{u}_s,\bm{v}_s)
=\rho(\bm{f}_s,\bm{v}_s)_{\Omega_s}
+(f_d,\psi_f)_{\Omega_d},\quad
\forall\,\underline{\bm{v}}\in\underline{\bm{W}},
\]
which implies that $\underline{\bm{u}}=(\bm{u}_s,\phi_m,\phi_f)\in\underline{\bm{W}}$ is a solution of \eqref{continuous:weak:formulation:2}.

Finally, we assume that there are two solutions $\underline{\bm{u}}^1,\underline{\bm{u}}^2\in\underline{\bm{W}}$ to \eqref{continuous:weak:formulation:2}. Then, there differences $\bm{\mathrm{e}}_s=\bm{u}_s^1-\bm{u}_s^2$, $\mathrm{e}_m=\phi_m^1-\phi_m^2$ and $\mathrm{e}_f=\phi_f^1-\phi_f^2$ satisfy that $\forall\,(\bm{v}_s,\psi_m,\psi_f)\in\bm{V}_s\times X_d\times X_d$,

\begin{equation}\label{error:equation}
\begin{split}
&2\rho\nu\left(\mathbb{D}(\bm{\mathrm{e}}_s),\mathbb{D}(\bm{v}_s)\right)_{\Omega_s}
+\frac{k_m}{\mu}\left(\nabla\mathrm{e}_m,\nabla\psi_m\right)_{\Omega_d}
+\frac{k_f}{\mu}\left(\nabla\mathrm{e}_f,\nabla\psi_f\right)_{\Omega_d}
+\left((\bm{\mathrm{e}}_s\cdot\nabla)\bm{u}_s^1,\bm{v}_s\right)_{\Omega_s}
+\left((\bm{u}_s^2\cdot\nabla)\bm{\mathrm{e}}_s,\bm{v}_s\right)_{\Omega_s}\\
&+\frac{\sigma k_m}{\mu}\left(\mathrm{e}_m-\mathrm{e}_f,\psi_m\right)_{\Omega_d}+\frac{\sigma k_m}{\mu}\left(\mathrm{e}_f-\mathrm{e}_m,\psi_f\right)_{\Omega_d}
+\left\langle\mathrm{e}_f,\bm{v}_s\cdot\bm{n}_s\right\rangle_\Gamma-
\left\langle\psi_f,\bm{\mathrm{e}}_s\cdot\bm{n}_s\right\rangle_\Gamma
+\sum_{i=1}^{N-1}\left\langle\bm{\mathrm{e}}_s\cdot\bm{\tau}_i,\bm{v}_s\cdot\bm{\tau}_i\right\rangle_\Gamma=0.
\end{split}
\end{equation}
Taking $\bm{v}_s=\bm{\mathrm{e}}_s$, $\psi_m=\mathrm{e}_m$ and $\psi_f=\mathrm{e}_f$ in \eqref{error:equation} and using another version of \eqref{trilinear:embedding}, which is

\begin{equation}\label{another:trilinear}
\left|\left((\bm{v}_s\cdot\nabla)\bm{w}_s,\bm{z}_s\right)_{\Omega_s}\right|
\leqslant \mathcal N\|\mathbb{D}(\bm{v}_s)\|_{0,\Omega_s}\|\mathbb{D}(\bm{w}_s)\|_{0,\Omega_s}\|\mathbb{D}(\bm{z}_s)\|_{0,\Omega_s},
\quad\forall\,\bm{v}_s,\bm{w}_s,\bm{z}_s\in\bm{X}_s,
\end{equation}
we obtain

\begin{equation}\label{error:equation:2}
  \left[2\rho\nu-\mathcal N(\|\mathbb{D}(\bm{u}_s^1)\|_{0,\Omega_s}+\|\mathbb{D}(\bm{u}_s^2)\|_{0,\Omega_s})\right]
  \|\mathbb{D}(\bm{\mathrm{e}}_s)\|_{0,\Omega_s}^2+\frac{k_m}{\mu}\|\nabla\mathrm{e}_m\|_{0,\Omega_d}^2+
  \frac{k_f}{\mu}\|\nabla\mathrm{e}_f\|_{0,\Omega_d}^2+\frac{\sigma k_m}{\mu}\|\mathrm{e}_m-\mathrm{e}_f\|_{0,\Omega_d}^2
  \leqslant 0.
\end{equation}
Hence, if we assume \eqref{assum:unique} holds, then \eqref{error:equation:2} shows the unique solution to \eqref{continuous:weak:formulation:2} based on a priori estimate \eqref{a:priori:estimate}.

\end{pf}

In order to prove the existence and uniqueness of the solution $(\underline{\bm{u}},p_s)\in\underline{\bm{Y}}\times Q_s$ to the model problem \eqref{continuous:weak:formulation}, we shall use the inf-sup condition \eqref{con:infsup} and the Babu\u{s}ka--Brezzi's theory \cite{Babuska1973179, Brezzi1974129, Girault1986, Temam1984}.

\begin{theorem}\label{theorem:unique:2}
  Under the assumption \eqref{assum:unique} of Theorem \ref{theorem:unique:1}, the model problem \eqref{continuous:weak:formulation} admits a unique solution $(\underline{\bm{u}},p_s)\in\underline{\bm{Y}}\times Q_s$ such that

  \begin{align}
  \label{estimate:velocity}
  &\rho\nu\|\mathbb{D}(\bm{u}_{s})\|_{0,\Omega_s}^2+\frac{2\sigma k_m}{\mu}\|\nabla\phi_{m}\|_{0,\Omega_d}^2
  +\frac{\sigma k_f}{\mu}\|\nabla\phi_{f}\|_{0,\Omega_d}^2\leqslant\mathcal C^2,\\
  \label{estimate:pressure}
  &\|p_s\|_{0,\Omega_s}\leqslant\beta_0^{-1}\left(2\rho^{1/2}\nu^{1/2}\mathcal C+\rho^{-1}\nu^{-1}\mathcal N\mathcal C^2+
  \rho\|\bm{f}_s\|_{-1,\Omega_s}+\|f_d\|_{-1,\Omega_d}\right),
  \end{align}
  where $\mathcal C^2$ is defined in Theorem \ref{theorem:priori:estimate}.
\end{theorem}

\begin{pf}
For the solution $\underline{\bm{u}}=(\bm{u}_s,\phi_m,\phi_f)\in \underline{\bm{W}}$ to \eqref{continuous:weak:formulation:2}, the following mapping:

\[
\underline{\bm{v}}=(\bm{v}_s,\psi_m,\psi_f)\in\underline{\bm{Y}}\;\mapsto\;
a(\underline{\bm{u}},\underline{\bm{v}})+b(\bm{u}_s;\bm{u}_s,\bm{v}_s)-\rho(\bm{f}_s,\bm{v}_s)_{\Omega_s}-(f_d,\psi_f)_{\Omega_d}
\]
defines an element $L(\underline{\bm{v}})$ of the dual space $\underline{\bm{Y}}'$, and furthermore, $L$ vanishes on $\underline{\bm{W}}$. As a result, the inf-sup condition \eqref{con:infsup} implies that there exists exactly one $p_s\in Q_s$ such that

\begin{equation}\label{unique:pressure}
L(\underline{\bm{v}})=d(\bm{v}_s,p_s),\quad\forall\,\underline{\bm{v}}=(\bm{v}_s,\psi_m,\psi_f)\in\underline{\bm{Y}}.
\end{equation}
Therefore the fact $\bm{u}_s\in\bm{X}_s$ and \eqref{unique:pressure} show that the model problem \eqref{continuous:weak:formulation} admits a unique solution $(\underline{\bm{u}},p_s)\in\underline{\bm{Y}}\times Q_s$. Finally, the result \eqref{estimate:pressure} is a straightforward application of the inf-sup condition \eqref{con:infsup} with the help of \eqref{another:trilinear} and \eqref{estimate:velocity}.

\end{pf}

%\section*{References}
\bibliographystyle{plain}
\bibliography{bibfile}

\begin{thebibliography}{}
\expandafter\ifx\csname url\endcsname\relax
  \def\url#1{\texttt{#1}}\fi
\expandafter\ifx\csname urlprefix\endcsname\relax\def\urlprefix{URL }\fi
\expandafter\ifx\csname href\endcsname\relax
  \def\href#1#2{#2} \def\path#1{#1}\fi

%\bibitem{achdou2000convergence}
%Y.~Achdou, J.~L. Guermond,
%Convergence Analysis of a Finite Element Projection/Lagrange--Galerkin Method for the Incompressible Navier--Stokes Equations,
%SIAM J. Numer. Anal. 37~(3) (2000) 799--826.
%\newblock \href {https://doi.org/10.1137/S0036142996313580}
%  {\path{doi:10.1137/S0036142996313580}}.

\bibitem{Adams2003}
R. Adams, J. Fournier,
Sobolev Spaces,
Acadamics Press, New York, 2003.

%\bibitem{Agmon2010}
%S. Agmon, Lectures on elliptic boundary value problems,
%AMS Providence, Ehode Island, 2010.
%\newblock \href {https://doi.org/10.1090/chel/369}
%  {\path{doi:10.1090/chel/369}}.

\bibitem{AlMahbub2019803}
M. A. Al Mahbub, X.-M. He, N. J. Nasu, C. Qiu, H. Zheng,
Coupled and decoupled stabilized mixed finite element methods for nonstationary dual--porosity--Stokes fluid flow model,
Int. J. Numer. Methods. Eng. 120 (6) (2019) 803--833.
\newblock \href {https://doi.org/10.1002/nme.6158}
  {\path{doi:10.1002/nme.6158}}.

\bibitem{AlMahbub2020112616}
M. A. Al Mahbub, F. Shi, N. J. Nasu, Y. Wang, H. Zheng,
Mixed stabilized finite element method for the stationary Stokes--dual--permeability fluid flow model,
Comput. Methods Appl. Mech. Engrg. 358 (2020) 112616.
\newblock \href {https://doi.org/10.1016/j.cma.2019.112616}
  {\path{doi:10.1016/j.cma.2019.112616}}.

%\bibitem{Arnold20021749}
%D. N. Arnold, F. Brezzi, B. Cockburn, L. D. Marini,
%Unified analysis of discontinuous Galerkin methods for elliptic problems,
%SIAM. J. Numer. Anal. 39 (5) (2002) 1749--1779.
%\newblock \href {https://doi.org/10.1137/S0036142901384162}
%  {\path{doi:10.1137/S0036142901384162}}.

\bibitem{Aubin1982}
T. Aubin,
Nonlinear Analysis on Manifolds. Monge-Amp\`{e}re Equations,
Grundlehren der mathematischen Wissenschaften 252,
Springer-Verlag New York, NY, 1982.
\newblock \href {https://doi.org/10.1007/978-1-4612-5734-9}
  {\path{doi:10.1007/978-1-4612-5734-9}}.

\bibitem{Babuska1973179}
I. Babu\u{s}ka,
The finite element method with Lagrangian multipliers,
Numer. Math. 20 (1973) 179--192.
\newblock \href {https://doi.org/10.1007/BF01436561}
  {\path{doi:10.1007/BF01436561}}.

\bibitem{Badea2010195}
L. Badea, M. Discacciati, A. Quarteroni,
Numerical analysis of the Navier-Stokes/Darcy coupling,
Numer. Math. 115 (2010) 195--227.
\newblock \href {https://doi.org/10.1007/s00211-009-0279-6}
  {\path{doi:10.1007/s00211-009-0279-6}}.

\bibitem{Beavers1967197}
G. Beavers, D. Joseph,
Boundary conditions at a naturally permeable wall,
Journal of Fluid Mechanics 30 (1) (1967) 197--207.
\newblock \href {https://doi.org/10.1017/S0022112067001375}
  {\path{doi:10.1017/S0022112067001375}}.

%\bibitem{becker2003a}
%R.~Becker, P.~Hansbo, R.~Stenberg, A finite element method for domain decomposition with non-matching grids,
%Mathematical Modelling and Numerical Analysis 37~(2) (2003) 209--225.
%\newblock \href {https://doi.org/10.1051/m2an:2003023}
%  {\path{doi:10.1051/m2an:2003023}}.

\bibitem{Bergh1976}
J. Bergh, J. L\"{o}fstr\"{o}m,
Interpolation Spaces: An Introduction,
Grundlehren der mathematischen Wissenschaften 223,
Springer-Verlag Berlin Heidelberg, Berlin, 1976.
\newblock \href {https://doi.org/10.1007/978-3-642-66451-9}
  {\path{doi:10.1007/978-3-642-66451-9}}.

%\bibitem{Bi2012425}
%C. Bi, M. Liu,
%A discontinuous finite volume element method for second--order elliptic problems,
%Numer. Methods Partial Differ. Equ. 28 (2) (2012) 425--440.
%\newblock \href {https://doi.org/10.1002/num.20626}
%  {\path{doi:10.1002/num.20626}}.

%\bibitem{Brenner2003306}
%S. C. Brenner,
%Poincar\'{e}-Friedrichs inequalities for piecewise $H^1$ functions,
%SIAM. J. Numer. Anal. 41 (1) (2003) 306--324.
%\newblock \href {https://doi.org/10.1137/S0036142902401311}
%  {\path{doi:10.1137/S0036142902401311}}.

%\bibitem{brenner2007mathematical}
%S. C. Brenner, L. R.~Scott,
%The mathematical theory of finite element methods, 3rd Edition,
%Texts in Applied Mathematics, Vol. 15, Springer, New York, NY, 2008.
%\newblock \href {https://doi.org/10.1007/978-0-387-75934-0}
%  {\path{doi:10.1007/978-0-387-75934-0}}.

\bibitem{Brezzi1974129}
F. Brezzi,
On the existence, uniqueness and approximation of saddle-point problems arising from Lagrangian multipliers,
RAIRO Anal. Numer. R2 (1974) 129--151.
\newblock \href {https://doi.org/10.1051/m2an/197408R201291}
  {\path{doi:10.1051/m2an/197408R201291}}.

%\bibitem{burman2007a}
%E.~Burman, P.~Hansbo, A unified stabilized method for Stokes' and Darcy's equations,
%Journal of Computational and Applied Mathematics 198~(1) (2007) 35--51.
%\newblock \href {https://doi.org/10.1016/j.cam.2005.11.022}
%  {\path{doi:10.1016/j.cam.2005.11.022}}.

\bibitem{Cai2009NUMERICAL}
M.~Cai, M.~Mu, J.~Xu, Numerical solution to a mixed {N}avier--{S}tokes/{D}arcy model by the two-grid approach,
SIAM J. Numer. Anal. 47~(5) (2009) 3325--3338.
\newblock \href {https://doi.org/10.1137/080721868}
  {\path{doi:10.1137/080721868}}.

\bibitem{Cao2021113128}
L. Cao, Y. He, J. Li, D. Yang,
Decoupled modified characteristic FEMs for fully evolutionary Navier-Stokes-Darcy model with the Beavers-Joseph interface condition,
Journal of Computational and Applied Mathematics 383 (2021) 113128.
\newblock \href {https://doi.org/10.1016/j.cam.2020.113128}
  {\path{doi:10.1016/j.cam.2020.113128}}.

\bibitem{Cao20132013}
Y. Cao, Y. Chu, X. He, M. Wei,
Decoupling the Stationary Navier-Stokes-Darcy System with the Beavers-Joseph-Saffman Interface Condition,
Abstract and Applied Analysis 2013 (2013) 136483.
\newblock \href {http://doi.org/10.1155/2013/136483}
  {\path{doi:10.1155/2013/136483}}.

\bibitem{Cao2011Robin}
Y.~Cao, M.~Gunzburger, X.~He, X.~Wang,
Robin--Robin domain decomposition methods for the steady-state {S}tokes--{D}arcy system with the {B}eavers--{J}oseph interface condition,
Numer. Math. 117~(4) (2011) 601--629.
\newblock \href {https://doi.org/10.1007/s00211-011-0361-8}
  {\path{doi:10.1007/s00211-011-0361-8}}.

\bibitem{Cao20104239}
Y. Cao, M. Gunzburger, X. Hu, F. Hua, X. Wang, W. Zhao,
Finite element approximations for Stokes--Darcy flow with Beavers--Joseph interface conditions,
SIAM. J. Numer. Anal. 47 (6) (2010) 4239--4256.
\newblock \href {https://doi.org/10.1137/080731542}
  {\path{doi:10.1137/080731542}}.

\bibitem{Cao20101}
Y. Cao, M. Gunzburger, F. Hua, X. Wang,
Coupled Stokes--Darcy model with Beavers--Joseph interface boundary condition,
Commun. Math. Sci. 8 (1) (2010) 1--25.
\newblock \href {https://doi.org/10.4310/CMS.2010.v8.n1.a2}
  {\path{doi:10.4310/CMS.2010.v8.n1.a2}}.

%\bibitem{chen2002error}
%Z.~Chen, R.~E. Ewing, Q.~Jiang, A.~M. Spagnuolo,
%Error analysis for characteristics-based methods for degenerate parabolic problems,
%SIAM J. Numer. Anal. 40~(4) (2002) 1491--1515.
%\newblock \href {https://doi.org/10.1137/S003614290037068X}
%  {\path{doi:10.1137/S003614290037068X}}.

\bibitem{Chen2012Efficient}
W.~Chen, M.~Gunzburger, S.~Dong, X.~Wang, Efficient and long-time accurate second-order methods for {S}tokes-{D}arcy system,
SIAM J. Numer. Anal. 51~(5) (2013) 2563--2584.
\newblock \href {https://doi.org/doi:10.1137/120897705}
  {\path{doi:doi:10.1137/120897705}}.

\bibitem{Chidyagwai20093806}
P. Chidyagwai, B. Rivi\`{e}re,
On the solution of the coupled Navier--Stokes and Darcy equations,
Comput. Methods Appl. Mech. Engrg. 198 (47) (2009) 3806--3820.
\newblock \href {https://doi.org/10.1016/j.cma.2009.08.012}
  {\path{doi:10.1016/j.cma.2009.08.012}}.

\bibitem{Cioranescu2016}
D. Cioranescu, V. Girault, K. R. Rajagopal,
Mechanics and Mathematics of Fluids of the Differential Type,
Advances in Mechanics and Mathematics, vol. 35. Springer International Publishing, Cham, Switzerland, 2016.
\newblock \href {https://doi.org/10.1007/978-3-319-39330-8}
  {\path{doi:10.1007/978-3-319-39330-8}}.

%\bibitem{Cockburn2002319}
%B. Cockburn, G. Kanschat, D. Sch\"{o}tzau, C. Schwab,
%Local Discontinuous Galerkin Methods for the Stokes System,
%SIAM J. Numer. Anal. 40 (1) (2002) 319--343.
%\newblock \href {https://doi.org/10.1137/S0036142900380121}
%  {\path{doi:10.1137/S0036142900380121}}.

%\bibitem{Cui2010824}
%M. Cui, X. Ye,
%Unified analysis of finite volume methods for the Stokes equations,
%SIAM J. Numer. Anal. 48 (3) (2010) 824--839.
%\newblock \href {https://doi.org/10.1137/090780985}
%  {\path{doi:10.1137/090780985}}.

\bibitem{d2011robust}
C. D'Angelo, P. Zunino, Robust numerical approximation of coupled {S}tokes' and
  {D}arcy's flows applied to vascular hemodynamics and biochemical transport,
  ESAIM: Mathematical Modelling and Numerical Analysis 45~(3) (2011) 447--476.
\newblock \href {https://doi.org/10.1051/m2an/2010062}
  {\path{doi:10.1051/m2an/2010062}}.

\bibitem{Discacciati2004phd}
M. Discacciati, Domain decomposition methods for the coupling of surface and groundwater flows,
Ph.D. thesis, Ecole Polytechnique F\'{e}d\'{e}rale de Lausanne, Switzerland, 2004.
\newblock \href {http://dx.doi.org/10.5075/epfl-thesis-3117}
  {\path{doi:10.5075/epfl-thesis-3117}}.

\bibitem{Discacciati2009315}
M. Discacciati, A. Quarteroni,
Navier-Stokes/Darcy Coupling: Modeling, Analysis, and Numerical Approximation,
Rev. Mat. Complut. 22 (2) (2009) 315--426.
\newblock \href {http://dx.doi.org/10.5209/rev_REMA.2009.v22.n2.16263}
  {\path{doi:10.5209/rev_REMA.2009.v22.n2.16263}}.

\bibitem{Discacciati2017571}
M. Discacciati, R. Oyarz\'{u}a,
A conforming mixed finite element method for the Navier--Stokes/Darcy coupled problem,
Numer. Math. 135 (2017) 571--606.
\newblock \href {http://dx.doi.org/10.1007/s00211-016-0811-4}
  {\path{doi:10.1007/s00211-016-0811-4}}.

\bibitem{Fang2020915}
J. Fang, P. Huang, Y. Qin,
A two-level finite element method for the steady-state Navier-Stokes/Darcy model,
J. Korean Math. Soc. 57 (4) (2020) 915--933.
\newblock \href {http://dx.doi.org/10.4134/JKMS.j190449}
  {\path{doi:10.4134/JKMS.j190449}}.

%\bibitem{fernandez2014explicit}
%M.~A. Fern\'{a}ndez, J.~Gerbeau, S.~Smaldone, Explicit coupling schemes for a
%  fluid-fluid interaction problem arising in hemodynamics, SIAM Journal on
%  Scientific Computing 36~(6) (2014) A2557–A2583.
%\newblock \href {https://doi.org/10.1137/130948653}
%  {\path{doi:10.1137/130948653}}.

\bibitem{Gao202001}
L. Gao, J. Li,
A decoupled stabilized finite element method for the dual-porosity-Navier–Stokes fluid flow model arising in shale oil,
Numer. Methods Partial Differential Eq. (2020) 1--18.
\newblock \href {https://doi.org/10.1002/num.22718}
  {\path{doi:10.1002/num.22718}}.

\bibitem{Girault1986}
V. Girault, P. A. Raviart,
Finite element methods for Navier--Stokes equations,
Springer-Verlag, Berlin, 1986.
\newblock \href {https://doi.org/10.1007/978-3-642-61623-5}
  {\path{doi:10.1007/978-3-642-61623-5}}.

\bibitem{Girault20092052}
V. Girault, B. Rivi\`{e}re,
DG Approximation of Coupled Navier--Stokes and Darcy Equations by Beaver--Joseph--Saffman Interface Condition,
SIAM J. Numer. Anal. 47~(3) (2009) 2052--2089.
\newblock \href {https://doi.org/10.1137/070686081}
  {\path{doi:10.1137/070686081}}.

\bibitem{Girault201328}
V. Girault, G. Kanschat, B. Rivi\`{e}re,
On the Coupling of Incompressible Stokes or Navier-Stokes and Darcy Flows Through Porous Media.
In: J. Ferreira, S. Barbeiro, G. Pena, M. Wheeler(eds) Modelling and Simulation in Fluid Dynamics in Porous Media.
Springer Proceedings in Mathematics \& Statistics, vol 28. Springer, New York, NY, 2013.
\newblock \href {https://doi.org/10.1007/978-1-4614-5055-9_1}
  {\path{doi:10.1007/978-1-4614-5055-9_1}}.

%\bibitem{Leopold2008}
%L.~Grinberg, G.~E. Karniadakis, A scalable domain decomposition method for
%  ultra-parallel arterial flow simulations, Communications in Computational
%  Physics~(4) (2008) 1151--1169.

\bibitem{hanspal2006numerical}
N.~Hanspal, A.~Waghode, V.~Nassehi, R.~Wakeman, Numerical analysis of coupled
  {S}tokes/{D}arcy flows in industrial filtrations, Transport in Porous Media
  64~(1) (2006) 73.
\newblock \href {https://doi.org/10.1007/s11242-005-1457-3}
  {\path{doi:10.1007/s11242-005-1457-3}}.

\bibitem{he2020an}
X. He, N.~Jiang, C.~Qiu, An artificial compressibility ensemble algorithm for a
  stochastic {S}tokes-{D}arcy model with random hydraulic conductivity and
  interface conditions, International Journal for Numerical Methods in
  Engineering 121~(4) (2020) 712--739.
\newblock \href {https://doi.org/10.1002/nme.6241}
  {\path{doi:10.1002/nme.6241}}.

\bibitem{He2015264}
X. He, J. Li, Y. Lin, J. Ming,
A domain decomposition method for the steady-state Navier-Stokes-Darcy model with Beavers-Joseph interface condition,
SIAM J. Sci. Comput. 37 (5) (2015), 264--290.
\newblock \href {https://doi.org/10.1137/140965776}
  {\path{doi:10.1137/140965776}}.

%\bibitem{heywood1990finite}
%J.~G. Heywood, R.~Rannacher,
%Finite-Element Approximation of the Nonstationary Navier–Stokes Problem. Part IV: Error Analysis for Second-Order Time Discretization,
%SIAM Journal on Numerical Analysis 27~(2) (1990) 353--384.
%\newblock \href {https://doi.org/10.1137/0727022}
%  {\path{doi:10.1137/0727022}}.

\bibitem{Hou2016710}
J. Hou, M. Qiu, X. He, C. Guo, M. Wei, B. Bai,
A dual--porosity--Stokes model and finite element method for coupling dual--porosity flow and free flow,
SIAM J. Sci. Comput. 38 (5) (2016) B710--B739.
\newblock \href {https://doi.org/10.1137/15M1044072}
  {\path{doi:10.1137/15M1044072}}.

\bibitem{Hou201947}
Y. Hou, S. Pei,
On the weak solutions to steady Navier-Stokes equations with mixed boundary conditions,
Mathematische Zeitschrift 291 (2019) 47--54.
\newblock \href {https://doi.org/10.1007/s00209-018-2072-7}
  {\path{doi:10.1007/s00209-018-2072-7}}.

%\bibitem{2017A}
%H. Jia, P. Shi, K. Li, H. Jia,
%A decoupling method with different subdomain time steps for the non-stationary Navier-Stokes/Darcy model,
%Journal of Computational Mathematics 35~(3) (2017) 319--345.
%\newblock \href {https://doi.org/10.4208/jcm.1606-m2015-0436}
%  {\path{doi:10.4208/jcm.1606-m2015-0436}}.

\bibitem{Jone1973231}
I. P. Jones,
Low Reynolds number flow past a porous spherical shell,
Mathematical Proceedings of the Cambridge Philosophical Society 73 (1) (1973) 231--238.
\newblock \href {https://doi.org/10.1017/S0305004100047642}
  {\path{doi:10.1017/S0305004100047642}}.

%\bibitem{kirby2003on}
%R.~C. Kirby, On the convergence of high resolution methods with multiple time
%  scales for hyperbolic conservation laws, Mathematics of Computation 72~(243)
%  (2003) 1239--1250.
%\newblock \href {https://doi.org/10.1090/S0025-5718-02-01469-2}
%  {\path{doi:10.1090/S0025-5718-02-01469-2}}.

%\bibitem{Kumar2015956}
%S. Kumar, R. Ruiz-Baier,
%Equal order discontinuous finite volume element methods for the Stokes Problem,
%J. Sci. Comput. 65 (2015), 956--978.
%\newblock \href {https://doi.org/10.1007/s10915-015-9993-7}
%  {\path{doi:10.1007/s10915-015-9993-7}}.

\bibitem{Li201692}
R. Li, J. Li, Z. Chen, Y. Gao,
A stabilized finite element method based on two local Gauss integrations for a coupled Stokes--Darcy problem,
Journal of Computational and Applied Mathematics 292 (15) (2016), 92--104.
\newblock \href {https://doi.org/10.1016/j.cam.2015.06.014}
  {\path{doi:10.1016/j.cam.2015.06.014}}.
%
%\bibitem{li2018a}
%Y.~Li, Y.~Hou, A second-order partitioned method with different subdomain time steps for the evolutionary Stokes-Darcy system,
%Math Meth Appl Sci. 41~(5) (2018) 2178--2208.
%\newblock \href {https://doi.org/10.1002/mma.4744}
%  {\path{doi:10.1002/mma.4744}}.
%
%\bibitem{Lions1972}
%J. L. Lions, E. Magenes,
%Non-homogeneous boundary value problems and applications, Vol. 1,
%Springer-Verlag, New York-Heidelberg, 1972.
%\newblock \href {https://doi.org/10.1007/978-3-642-65161-8}
%  {\path{doi:10.1007/978-3-642-65161-8}}.

\bibitem{Mikelic20001111}
A. Mikeli\'{c}, W. J\"{a}ger,
On the interface boundary condition of Beavers, Joseph, and Saffman,
SIAM Journal on Applied Mathematics 60 (4) (2000) 1111--1127.
\newblock \href {https://doi.org/10.1137/S003613999833678X}
  {\path{doi:10.1137/S003613999833678X}}.

\bibitem{Mu2010Decoupled}
M.~Mu, X.~Zhu, Decoupled schemes for a non-stationary mixed {S}tokes-{D}arcy model,
Mathematics of Computation 79~(270) (2010) 707--731.
\newblock \href {https://doi.org/10.1090/S0025-5718-09-02302-3}
  {\path{doi:10.1090/S0025-5718-09-02302-3}}.

%\bibitem{Nitsche1971}
%J.~Nitsche, \"{U}ber ein Variationsprinzip zur L\"{o}sung von Dirichlet-Problemen bei Verwendung von Teilr\"{a}umen, die keinen Randbedingungen unterworfen sind,
%Abhandlungen aus dem Mathematischen Seminar der Universit\"{a}t Hamburg 36 (1971) 9--15.
%\newblock \href {https://doi.org/10.1007/BF02995904}
%  {\path{doi:10.1007/BF02995904}}.

\bibitem{nassehi1998modelling}
V.~Nassehi, Modelling of combined {N}avier-{S}tokes and {D}arcy flows in
  crossflow membrane filtration, Chemical Engineering Science 53~(6) (1998)
  1253--1265.
\newblock \href {https://doi.org/10.1016/S0009-2509(97)00443-0}
  {\path{doi:10.1016/S0009-2509(97)00443-0}}.

\bibitem{Qiu2020109400}
C. Qiu, X. He, J. Li, Y. Lin,
A domain decomposition method for the time-dependent Navier-Stokes-Darcy model with Beavers-Joseph interface condition and defective boundary condition,
Journal of Computational Physics 411 (2020) 109400.
\newblock \href {https://doi.org/10.1016/j.jcp.2020.109400}
  {\path{doi:10.1016/j.jcp.2020.109400}}.

\bibitem{Riviere2008}
B. Rivi\`{e}re,
Discontinuous Galerkin methods for solving elliptic and parabolic equations.
Theory and implementation, Front. Appl. Math. 35, SIAM, Philadelphia, 2008.
\newblock \href {http://dx.doi.org/10.1137/1.9780898717440}
  {\path{doi:10.1137/1.9780898717440}}.

\bibitem{Roos2008}
H. G. Roos, M. Stynes, L. Tobiska,
Robust numerical methods for singularly perturbed differential equations.
Convection--Diffusion--Reaction and Flow Problems, 2nd ed., Springer Ser. Comput. Math. 24, Springer, Berlin, 2008.
\newblock \href {https://doi.org/10.1007/978-3-540-34467-4}
  {\path{doi:10.1007/978-3-540-34467-4}}.

%\bibitem{rybak2014a}
%I.~Rybak, J.~Magiera, A multiple-time-step technique for coupled free flow and
%  porous medium systems, Journal of Computational Physics 272 (2014) 327--342.
%\newblock \href {https://doi.org/10.1016/j.jcp.2014.04.036}
%  {\path{doi:10.1016/j.jcp.2014.04.036}}.

\bibitem{Saffman197193}
P. G. Saffman,
On the boundary condition at the surface of a porous medium,
Studies in Applied Mathematics 50 (2) (1971) 93--101.
\newblock \href {https://doi.org/10.1002/sapm197150293}
  {\path{doi:10.1002/sapm197150293}}.

\bibitem{Scott1990483}
L. Scott, S. Zhang,
Finite element interpolation of nonsmooth functions satisfying boundary conditions,
Math. Comp. 54 (190) (1990) 483--493.
\newblock \href {https://doi.org/10.1090/S0025-5718-1990-1011446-7}
  {\path{doi:10.1090/S0025-5718-1990-1011446-7}}.

\bibitem{New1983}
K. Serra, A. C. Reynolds, R. Raghavan,
New pressure transient analysis methods for naturally fractured reservoirs,
Journal of Petroleum Technology 35~(12) (1983) 2271--2283.
\newblock \href {https://doi.org/10.2118/10780-PA}
  {\path{doi:10.2118/10780-PA}}.

\bibitem{Shan2019389}
L. Shan, J. Hou, W. Yan, J. Chen,
Partitioned time stepping method for a dual--porosity--Stokes model,
J. Sci. Comput. 79 (1) (2019) 389--413.
\newblock \href {https://doi.org/10.1007/s10915-018-0879-3}
  {\path{doi:10.1007/s10915-018-0879-3}}.

\bibitem{Shan2013813}
L. Shan, H. Zheng,
Partitioned time stepping method for fully evolutionary Stokes--Darcy flow with Beavers--Joseph interface conditions,
SIAM J. Numer. Anal. 51 (2) (2013) 813--839.
\newblock \href {https://doi.org/10.1137/110828095}
  {\path{doi:10.1137/110828095}}.
%
%\bibitem{shan2013a}
%L.~Shan, H.~Zheng, W.~Layton,
%A decoupling method with different subdomain time steps for the nonstationary Stokes--Darcy model,
%Numerical Methods for Partial Differential Equations 29~(2) (2013) 549--583.
%\newblock \href {https://doi.org/10.1002/num.21720}
%  {\path{doi:10.1002/num.21720}}.
%
%\bibitem{Si2014Decoupled}
%Z.~Si, Y.~Wang, S.~Li, Decoupled modified characteristics finite element method for the time dependent {N}avier--{S}tokes/{D}arcy problem,
%Mathematical Methods in the Applied Sciences 37~(9) (2014) 1392--1404.
%\newblock \href {https://doi.org/10.1002/mma.2901}
%  {\path{doi:10.1002/mma.2901}}.
%
%\bibitem{Si2016Unconditional}
%Z.~Si, J.~Wang, W.~Sun, Unconditional stability and error estimates of modified
%  characteristics {FEM}s for the {N}avier--{S}tokes equations, Numerische
%  Mathematik 134~(1) (2016) 139--161.
%\newblock \href {https://doi.org/10.1007/s00211-015-0767-9}
%  {\path{doi:10.1007/s00211-015-0767-9}}.

\bibitem{Temam1984}
R. Temam,
Navier--Stokes Equations: Theory and Numerical Analysis, 3rd Edition,
North-Holland, Amsterdam, 1984.
\newblock \href {https://doi.org/10.1090/chel/343}
  {\path{doi:10.1090/chel/343}}.

\bibitem{The1963}
J. E. Warren, P. J. Root, The behavior of naturally fractured reservoirs, Society of
  Petroleum Engineers Journal 3~(3) (1963) 245--255.
\newblock \href {https://doi.org/10.2118/426-PA}
  {\path{doi:10.2118/426-PA}}.

%\bibitem{Ye20041062}
%X. Ye,
%A new discontinuous finite volume method for elliptic problems,
%SIAM J. Numer. Anal. 42 (3) (2004) 1062--1072.
%\newblock \href {https://doi.org/10.1137/S0036142902417042}
%  {\path{doi:10.1137/S0036142902417042}}.
%
%\bibitem{Ye2006183}
%X. Ye,
%A discontinuous finite volume method for the Stokes problems,
%SIAM J. Numer. Anal. 44 (1) (2006) 183--198.
%\newblock \href {https://doi.org/10.1137/040616759}
%  {\path{doi:10.1137/040616759}}.
%
%\bibitem{Yin2012761242}
%Z. Yin, Z. Jiang, Q. Xu,
%A Discontinuous Finite Volume Method for the Darcy--Stokes Equations,
%Journal of Applied Mathematics 2012 (2012) 761242.
%\newblock \href {https://doi.org/10.1155/2012/761242}
%  {\path{doi:10.1155/2012/761242}}.

\bibitem{zhao2016two-grid}
J.~Zhao, T.~Zhang, Two-grid finite element methods for the steady {N}avier-{S}tokes/{D}arcy model,
East Asian Journal on Applied Mathematics 6~(1) (2016) 60--79.
\newblock \href {https://doi.org/10.4208/eajam.080215.111215a}
  {\path{doi:10.4208/eajam.080215.111215a}}.

\bibitem{Zuo20151009}
L. Zuo, Y. Hou,
Numerical analysis for the mixed Navier--Stokes and Darcy Problem with the Beavers--Joseph interface condition,
Numer. Methods Partial Differential Eq. 31 (2015) 1009--1030.
\newblock \href {https://doi.org/10.1002/num.21933}
  {\path{doi:10.1002/num.21933}}.

\end{thebibliography}
\end{document}